\documentclass[reqno]{amsart}
\usepackage{amsmath}
\usepackage{amsthm}
\usepackage{amssymb}
\usepackage{graphics}

\DeclareMathOperator{\supp}{supp}

\usepackage{pst-all}
\usepackage{caption2}

\usepackage[latin1]{inputenc} %%% Para poder escribir áéíóúñ directamente!

\theoremstyle{definition}
\newtheorem{teorema}{Theorem}
\newtheorem{coro}{Corollary}

\newtheorem{propo}{Proposition}
\newtheorem{lema}{Lemma}[section]

\newtheorem*{notaciones}{Notation and definitions}
\newtheorem*{nota}{Remark}

\newtheorem{afir}{Claim}[section]

\newcommand{\R}{\mathbb R}

\newcommand{\Z}{\mathbb Z}
\newcommand{\N}{\mathbb N}
\newcommand{\g}{\gamma}

\newcommand{\s}{\sigma}
\renewcommand{\b}{\beta}

\renewcommand{\d}{\delta}
\renewcommand{\l}{\lambda}
\renewcommand{\t}{\theta}
\renewcommand{\a}{\alpha}
\renewcommand{\o}{\omega}

\newcommand{\q}{\hat{q}}

\renewcommand{\t}{\theta}

\newcommand{\suma}{\displaystyle\sum}
\newcommand{\union}{\displaystyle\bigcup}
\newcommand{\maxim}{\displaystyle\max}

\newcommand{\norm}[1]{\| #1 \|}

\renewcommand{\qed}{\hskip1em \Box}

\newcommand{\um}{{1\over 2}}

\newcommand{\tupla}{(\xi , \tau , \phi)}
\newcommand{\dominio}{\Gae^3}

\newcommand{\e}{|\eta|}

\newcommand{\xim}{|\xi|}
\newcommand{\taum}{|\tau|}
\newcommand{\phim}{|\phi|}

\newcommand{\esca}{\xi\cdot (\eta-\xi)}

\newcommand{\expdos}{\q(\eta-\tau)}

\newcommand{\etau}{|\eta-\tau|}
\newcommand{\etam}{{\eta\over 2}}
\newcommand{\metam}{{\e\over 2}}
\newcommand{\xitau}{|\xi-\tau|}
\newcommand{\metac}{{\e\over 100}}

\newcommand{\qephitau}{\q(\tau'+\phi'-\eta)}

\newcommand{\Gae}{\Gamma (\eta )}
\newcommand{\Lax}{\Lambda (\xi)}

\newcommand{\Gainf}{\Gamma_{\infty} (\eta )}

\newcommand{\homoa}{\dot{W}^{\a,2}}
\newcommand{\homob}{\dot{W}^{\b,2}}
\newcommand{\homog}{\dot{W}^{\g,2}}

\newcommand{\nohomoa}{W^{\a,2}}

\newcommand{\homobmem}{\dot{W}^{\b-{1\over 2}+\epsi\, ,\, 2}}

\newcommand{\homomenos}{\dot{W}^{-{1\over 2},2}}

\renewcommand{\qed}{\hfill$\Box$}

\newcommand{\proofof}[1]{ \noindent {\bf Proof of #1. }}

\newcommand{\rteta}{R_{\t , k}}

\newcommand{\f}[2]{\frac{#1}{#2}}
\newcommand{\Rn}{\mathbb{R}^{n}}
\newcommand{\ra}{\rightarrow}

\newcommand{\Rtres}{\mathbb{R}^{3}}

\newcommand{\reg}{C^{\infty}(\Rn)}

\newcommand{\wt}{\widetilde}
\newcommand{\wh}{\widehat}

\newcommand{\etauphi}{|\eta-\tau-\phi|}

\newcommand{\exiphi}{|\eta-\xi-\phi|}

\newcommand{\F}{\mathcal{F}}

\newcommand{\Gad}{\Gamma_{\d}(\eta)}

\newcommand{\fourqq}{\q(\xi)\q(\eta-\tau)\q(\tau-\phi')\q(\phi'-\xi)}
\newcommand{\fourq}{\q(\xi)\q(\eta-\tau)\q(\tau+\phi-\eta)\q(\eta-\phi-\xi)}
\newcommand{\measuree}{d\s(\xi)d\s(\tau)d\s(\phi')}
\newcommand{\measure}{d\s(\xi)d\s(\tau)d\s(\phi)}
\newcommand{\measureee}{d\s(\xi')d\s(\tau')d\s(\phi')}

\newcommand{\cte}{{\sqrt{2}\over 2}+{1\over 100}}

\allowdisplaybreaks

\newcommand {\epsi} {\varepsilon}
\newcommand {\ga} {{\gamma}}

\numberwithin{equation}{section}

\begin{document}

\title[Reconstruction of the singularities of a potential from backscattering data]
{Reconstruction of the singularities of a potential from
backscattering data in 2D and 3D}

\thanks{Both authors were supported   by Spanish Grant
MTM2005-07652-C02-01.}

\author[Juan Manuel Reyes and Alberto Ruiz]{Juan Manuel Reyes$^1$ and Alberto Ruiz$^2$}

\maketitle

\begin{center}
{\small $^1$Departamento de Matem\'aticas y Computaci\'on\\ Universidad de La Rioja\\
 26004 Logro\~no, Spain\\
E-mail: juanmanuel.reyes@unirioja.es }
\end{center}

\begin{center}
{\small $^2$Departamento de Matem\'aticas\\ Universidad Aut\'onoma de Madrid\\
 28049 Madrid, Spain\\
E-mail: alberto.ruiz@uam.es }
\end{center}

\begin{abstract}
We prove  that the  singularities of  a potential in  the two and
three dimensional  Schr\"odinger equation  are  the  same as  the
singularities of  the Born approximation (Diffraction Tomography),
obtained from backscattering inverse data,  with an accuracy of
$1/2^-$  derivative in the  scale  of $L^2$-based Sobolev  spaces.
This improves  previous  results, see \cite{RV} and \cite{OPS},
removing  several constrains on  the  a priori  regularity  of the
potential. The improvement is based on the study of the smoothing
properties of the quartic term in the Neumann-Born expansion of
the scattering amplitude in 3D, together with a Leibniz formula
for multiple scattering valid  in  any dimension.
\end{abstract}

$ $\\ \\
\section{Introduction}
The inverse scattering problem for Schr\"odinger  potentials  deals with the  uniqueness, reconstruction and  stability of the potential $q $  in the  Hamiltonian
$H=\Delta  +q$  from  the  far field pattern  of the  generalized eigenfunctions or  scattering  solutions.  These  are the  unique  solutions  of  the asymptotic boundary  value problem
\begin{equation}\label{scattsol}
\begin{cases}
(\Delta+ q+k^2) u = 0\\
u = e^{ikx\cdot\theta}+  u_{out},
\end{cases}
\end{equation}
where the function  $u_{out}$ satisfies  the outgoing Sommerfeld
radiation condition, which means, for compactly supported
potential $q$,  that $u$  has  asymptotics  as $|x|\to \infty$
\begin{equation}\label{assymp}
u(x, \theta, k)= e^{ikx.\theta} + C
|x|^{\frac{1-n}{2}}k^{\frac{n-3}{2}}e^{ik|x|}A( \theta', \theta,k)
+o(|x|^{\frac{1-n}{2}}),
\end{equation}
where $\theta'=x/|x|$.

The function $A( \theta',\theta,k)$, $k > 0$,
$\theta,\theta'$ in the
 unit
 sphere $S^{n-1} $, is known as  scattering amplitude or far  field pattern.

 The  inverse problem for whole  data is formally overdeterminate,   as   one easily can see   by  counting variables.  For
 this reason, to  avoid   redundancies,  some  kinds of  partial  data   are selected for the inverse problems.  The  selection of
  these data is motivated by numerical experience and applications. The most celebrated sets of partial data are the 
   following:
 \begin{itemize}
 \item{Fixed energy  data.}  We assume  as  data  the  values  $A( \theta',\theta,k)$
 for  fixed $k$  and  free  $\theta,\theta'\in \mathbb S^{n-1}$. Uniqueness  of  the inverse problem  in
 this  case was studied  by \cite{Na}, \cite{No}, \cite{R}, \cite{SU2}.  The  approach to  this problem  is  related  to
 the Calder\'on- Sylvester-Uhlmann  complex  exponential solutions, used  in  the Electrical Impedance  Tomography
  inverse problem. The  stability  happens to  be  very  weak.
 \item{Fixed angle  data.} The knowledge  of  $A( \theta',\theta,k)$
 for  fixed $\theta $, free  $ \theta'\in \mathbb S^{n-1}$  and free $k>0$ is  assumed.
 The uniqueness of  the inverse problem is open  and only  generic  and  local  uniqueness is proved under a
 priori regularity assumptions on the potential, see \cite{S}.
 \item{Backscattering  data.}  One  assumes $A( -\theta,\theta,k)$ for free $ \theta\in \mathbb S^{n-1}$ and free
 $k>0$. The uniqueness of the inverse problem is   only proved generically and for small potentials, see \cite{ER1}, see also
\cite{P}, \cite{S}, \cite{L}.

\end{itemize}

 In practical applications the actual potential is substituted by the so called Born approximation of the scattering amplitude. The procedures to imaging
 the Born approximation from the scattering data are known as Diffraction Tomography.

 The  different Born approximations  are obtained, in the frequency domain,  from  the formula
 $$\hat q_{approx}(\xi)= A(\omega,\theta,k), $$
 where $\xi$ is given  by the redundant relation
 $$\xi=k(\omega-\theta).$$
If $\theta$ is fixed (fixed angle data), we use the change of
variable
 $\xi= \Phi_\theta(k,\omega)= k(\omega-\theta)$ to define  the Born approximation
 $\hat q_\theta(\xi) =A(\omega,\theta,k)$. Notice that this change of variable becomes singular   on the hyperplane $\xi\cdot \theta=0$.

For backscattering  data the Born approximation is given in the frequency domain by the polar coordinates
 \begin{equation}\label{bsborn}
 \hat q_B(-2k\theta)= A(-\theta,\theta,k).   
  \end{equation}
This fact makes backscattering  data more natural and simpler  than fixed angle data
for  diffraction tomography.

The use of the Born approximation is  not,  in  general, justified
on a mathematical  basis: one  would  like  to  know how much
information on the actual potential $q$ is contained in the Born
approximation. This problem has been treated by several authors.
In  full data case see \cite{PSe}, \cite{PS}, \cite{PSS} and
\cite{BFRV}. For fixed angle data  and backscattering data,  both of which are formally well determinate,  the justification of diffraction tomography  was studied in    \cite{R1}    (fixed angle) and  \cite{GU}, \cite{OPS}, \cite{RV}, \cite{Re1}    (backscattering). We would like to remark that each of the  last two types of data require the analysis of special multilinear operators which are not related.

In this  work we study the  case of backscattering
data in dimension two and three, we continue  and complete the
research  of \cite{OPS},  \cite{RV}  and \cite{Re1}, by removing
some constrains in their results.    We prove that  the diffraction
tomography is a migration scheme, see \cite{B}, within an accuracy
of at least 1/2. This  is to say that the most singular parts of
the  actual potential  can be reconstructed from   the  Born approximation up  to a
certain order (the accuracy of the migration). The determination  of
this accuracy  is very important to design numerical methods,
adapted to the  spaces in which one expects to obtain the
information on the actual potential from real scattering data. We prove
 \begin{teorema}
\label{theorem3}{\it Assume that $n\in\{2,3\}$, $q$ is a compactly
supported function in $W^{\a\, ,\, 2}(\Rn)$ and $\a\geq 0$. Then
$q-q_{B}\in W^{\b\, ,\, 2}(\Rn)+C^{\infty}(\Rn)\,$, for any
$\b\in\R$ such that $0\leq \b<\a+{1\over 2}\,$. }
\end{teorema}

In   Theorem \ref{theorem3} the regularity is measured in the scale of $L^2$-based Sobolev  spaces. The optimality of this accuracy in this  scale of spaces is, so far, an open and interesting question.

The procedure to justify  the migration scheme is  to study  the
smoothing properties  of  the  multi-linear terms in  the
Neumann-Born  expansion  of  the  scattering  amplitude (multiple
scattering).

Physical  evidence suggests that  multiple scattering  is strong
in the  case  of backscattering  data. The control of  double and
triple scattering in 3D, within an accuracy of 1/2,  was obtained
in \cite{RV}, but  their  result together with the
 general estimates for multiple scattering do not suffice to assure that, for  a potential $q$ a priori in
 the Sobolev space $W^{\alpha ,2}$, the error $q-q_{B}$ is in $W^{\beta,2}$ for  any  $\beta<\a + 1/2$; the
 restriction $0\le \alpha <3/4$ is needed. In  the    range
  $  \alpha \geq 3/4$ known estimates of quadruple  scattering  became worse  than those of double  or triple  scattering.  The study of  the accuracy of the Born approximation requires, then,  to improve the estimates of  the quartic term in the series.
  We accomplish this in  the present work.  We also  extend the results, which previously were only studied for $\alpha<3/2$ in 3D and for $\alpha<1$ in 2D, to  any $\alpha\geq 0$ by   using  a Leibniz' type formula for the derivatives of
  multiple scattering terms (see \S C.1 in \cite{Re2}).

 In  dimension three, we only are  able to prove that the  errors due  to  double, triple  and quadruple scattering are   a half of a derivative better than the actual potential, as opposite  to  the  2D case where the regularity increases with the order.

  Result   from  \cite{RV},  \cite{Re1}  together  with  Corollary \ref{lcontinua} allow us  to state  the following  result concerning  reconstruction  of classical discontinuities from backscattering in 2D:

 \begin {teorema}{\it Let $q$  compactly supported in $ W^{\alpha, 2}({\mathbb R}^2)$, where $\alpha>0 $. Then  $q-q_B $
 is  a continuous function.}
 \end{teorema}

 In fact, it was
proved  (Theorem 2 in \cite{RV}), that, for such a $q$, the
quadratic term is a continuous function. H\"older continuity of
the cubic term is obtained since it  is in $ W^{\beta,2}$ for all
$\beta<\alpha +1$, see Theorem 1 in \cite{Re1}. The remainder is controled by Corollary \ref{lcontinua}.

In the three dimensional case, it follows  from Theorem 1  that  the  whole non continuous part of  the  actual potential
can be  reconstructed from the Born  approximation,  assuming  
a priori that  $q$ is  in the  Sobolev  space $W^{\alpha,2}$ for  some
$\alpha>1$. Notice $q$ might have  some discontinuities   if  $\alpha$ is between $ 1$ and the  3D Sobolev exponent $ 3/2$:

\begin{teorema}{\it
Let $q$  compactly supported in $W^{\alpha, 2}({\mathbb R}^3)$,
where $\alpha>1$. Then  $q-q_B $ is  a continuous  function.}
\end{teorema}

 From the  previous work  \cite{RV} it follows  that  in 3D the discontinuities in the case of a piecewise regular potential can be reconstructed from the Born approximation (the result is not stated in \cite{RV} but it is similar to Corollary 0.1 in \cite{OPS} in the  2D case).  By using  the evolution equation  the reconstruction of conormal
singularities was achived in \cite{GU}. On one hand Theorem 3, as far as  we know,   is the  first   result of reconstruction  of
discontinuities in 3D, without assuming  special  structure of the
singular set  but, on the other hand,
one  expects that  $q\in W^{\alpha,
2}$, for  any $\alpha>1/2$ suffices for the reconstruction of
discontinuities. So  far, this improvement has not been
achieved. We know from  Corollary \ref{lcontinua}, that the high frequency
Neumann-Born series for $j\geq 5$ converges to a H\"older
continuous function for $\a>\um$.

An important feature of Theorem \ref{theorem3} is the fact that, regardless of the a priori
regularity assumptions on the potential,  the accuracy of the migration scheme is at least $1/2$. This independency is important to construct any
 recurrence  scheme, in order to obtain further information on the actual potential from  scattering
 data. In the case of fixed angle data, one can define a modified
 Born approximation by inserting the error $q-q_B$ in the
 quadratic form, see \cite{R1}. This increases the known accuracy
 for rough potentials $q$, but an inconvenient to iterate the procedure is the dependency
 on $\a$ of the accuracy.  
 
Finally,  we remark that in the higher dimensional  case,  the  order of accuracy is  an open question. We believe  that 1/2 also applies, but
 the technical complexity of  our approach makes it  necessary   to look for   a new  point  of  view on  the problem. The  treatment of the 3D problem   due to
 Lagergren  \cite{L}, \cite{M},  based upon a time  dependent expansion  of  the  backscattering  operator,  also  requires  a very technical   treatment  of  its multilinear  term. See also \cite{BeM1}	 and \cite{BeM2}.

\begin{notaciones}
We will write $\F f$ or $\hat{f}$ to denote the Fourier transform
of $f$. $\F^{-1}$ denotes the inverse Fourier transform. The
letter $M$ denotes the Hardy-Littlewood
 maximal operator. We denote the two-dimensional
 Hausdorff measure in $\Rtres$ by $\s$. The expression $|x|\sim 2^{-j}\e$ refers to $2^{-j-1}\e <|x|\leq 2^{-j+1}\e$, for
 $j\in \mathbb{Z}$, $x,\eta\in\R^3$. We will use the homogeneous and non homogeneous Hilbertian Sobolev
spaces. With $\a\in\R$, we denote
\begin{align*}
& \dot{W}^{\a , 2}(\Rn):=\{f\in \mathcal{S}'(\Rn) : |\cdot|^{\a}\F
f(\cdot)\in L^2(\Rn)\},\\
& W^{\a , 2}(\Rn):=\{f\in \mathcal{S}'(\Rn) : (1+|\cdot|^2)^{\a /
2}\F f (\cdot)\in L^2(\Rn)\}.
\end{align*}

 Let
$\eta,\, \xi\in\Rtres\setminus\{0\}\,$. We write
\begin{equation}
\label{notation1}\Gae :=\left\{ x \in\Rtres : \left| x -{\eta\over
2}\right|={\e\over 2}\right\}\, ,
\end{equation}
refering to the sphere centered at ${\eta\over 2}$ and radius
${\e\over 2}\,$ and
\begin{equation}
\label{notation2}\Lambda(\xi):=\left\{ x \in\Rtres : \xi\cdot ( x
-\xi )=0 \right\}
\end{equation}
denotes the plane orthogonal to $\xi$ that contains the point
$\xi$. We denote $\dominio:=\Gae\times\Gae\times\Gae$ and
$\Gae^2:=\Gae\times\Gae$.

Let $\widehat{F}(\eta)$ given  by  the integral on a manifold
$A(\eta)$ of some function. Since our proofs are based upon a
decomposition of $A(\eta)$ in several subdomains $D(\eta) \subset
A(\eta)$, we will denote by $ \widehat{F}_D(\eta)$ the same
expression when  we restrict the integration  to the subdomain
$D(\eta)$.

%Let $F$ be defined by $\widehat F(\eta)
%=\int_{A(\eta)}f(\eta,\xi)d\sigma(\xi)$. Since our proofs are
%based upon a  decomposition of  $A(\eta)$  in several subdomains
%$D(\eta)\subset A(\eta)$, we will denote $\widehat F_D(\eta)
%=\int_{D(\eta)}f(\eta,\xi)d\sigma(\xi)$.

 The outgoing resolvent operator for the Laplacian is defined, in
terms of the Fourier transform, by
\[
\widehat{R_{k} (f)}(\xi)=  (-\xim^2 + k^2 + i0 )^{-1}
\widehat{f}(\xi)\, .
\]
We define the operator $Q_j$ in the following way
\begin{equation}
\label{term}\widehat{Q_j (q)}(\xi) = \int_{\Rn} e^{ik\t\,\cdot\,
y} (q R_k )^{j-1} (q(\cdot ) e^{ik\t\cdot (\cdot)})(y)dy\, ,
\end{equation}
where  $k=\xim / 2$, $\t= -\xi / \xim$. With these expressions for $k$
and $\t$, we define the multi-linear form $Q_j(f_1 , \,\dots\, ,
f_j)$ in the FT side as
\begin{align}
\notag&\mathcal{F}\left( Q_j (f_1 , \dots , f_{j}) \right)
(\xi):=\int_{\Rn}e^{ik\t\cdot
y}(f_1R_k)(f_2R_k)\,\dots\,(f_{j-1}R_k)(f_j(\cdot)e^{ik\t\cdot
(\cdot)})(y)dy .
\end{align}
We denote the high frequency version
\begin{equation}\label{qtilde}
\widetilde{Q}_{j}(q):=\mathcal{F}^{-1}(\chi^{\ast}(|\cdot|/2)\widehat{Q_{j}(q)}(\cdot))\,
,
\end{equation}
   where $Q_j(q)$ is defined in \eqref{term} and $\chi^{\ast}\in C^{\infty}(\R)$ with
$\chi^{\ast}(t)=1$ if $ t\geq 2C_0$, $\chi^{\ast}(t)=0$ if $t <
C_0$, for a certain constant $C_0>0$ to be chosen (see section 4).
Notice that the cutoff near the origin allows us to reduce the
estimates of Sobolev norms to the estimation of homogeneous
Sobolev  norms.

We also write $\wt{Q}_j(f_1 , \dots , f_j)=\mathcal{F}^{-1}\left(
\chi^{\ast} \mathcal{F}\left( Q_j(f_1 , \dots , f_j ) \right)\,
\right)\,$. We will admit the abuse of notation $Q_j(q)=Q_j(q,\,
\dots \, , q)$ and $\wt{Q}_j(q)=\wt{Q}_j(q,\, \dots \, , q)$.

The permutation group of order $k$ is denoted by $S_k$. For
multi-indexes $\b$ and $\ga$ in $\mathbb{N}^{n}\,$, we use the
standard definitions  of $\b!$, $|\b|$ and  $\b\leq \ga$.

We use the letter $C$ to denote any constant that can be
explicitly computed in terms of known quantities. The exact value
denoted by $C$ may change from line to line in a given
computation.

\end{notaciones}

\section{Preliminaries and results}

We obtain the  so called Lippmann-Schwinger integral equation by
applying the outgoing resolvent to \eqref{scattsol}
\begin{equation}
\label{lippschwinger}u(x , \t , k )= e^{ik x\cdot \,\t } +
R_{k}(q(\cdot )u(\cdot \, , \t , k ))(x)\, .
\end{equation}
The key operator in the above integral equation is
\[
T_{k}(f)(x) = R_{k}(q(\cdot )f(\cdot ))(x)\, .
\]
There are several a priori estimates for $R_{k}$ that allow to
prove existence and uniqueness of Lippmann-Schwinger integral
equation. Usually,  Fredholm theory applies and everything follows
from compactness arguments, Rellich uniqueness theorem and unique
continuation principles, in the case of real valued potentials.
The solution can be obtained in several situations (these cases do
not require $q$ to be  real) by perturbation arguments, assuming
that the energy is sufficiently large, $k>k_{0}\geq 0\,$, where
$k_{0}$ depends on some a priori bound of the potential $q\,$. As
an example we may consider compactly supported $q\in L^r (\Rn)$
for some $r>{n\over 2}\,$. In this case, which is  the one considered in this work, the resolvent operator
$R_k $ is bounded from $L^p (\Rn)$ to $L^{p' }(\Rn )$ with norm
decaying to $0$ as $k\ra \infty $ when ${1\over p} - {1\over p' }
= {1\over r}\,$, see \cite{A}, \cite{KRS} and see also \cite{R1}.
This together with Hölder inequality proves that for big $k$ the
operator $T_{k}$ is a contraction in $L^p $ and then existence and
uniqueness of solution of \eqref{lippschwinger} easily follow and $u$  can be expressed as  a convergent Neumann-Born series.

Once the scattering solution is obtained we may prove that the far
field pattern can be expressed as
\begin{equation}
\label{farfield}A(\t' , \t , k )=\int_{\Rn} e^{-ik\t' \cdot \, y}
q(y) u(y , \t , k) dy\, ,
\end{equation}
see \cite{ER1} where this  is  used  as  a definition  for
non compactly supported potentials.

By inserting the series   $u$  in \eqref{farfield} one obtains the Neumann-Born series of the scattering amplitude for $k$ large enough  (high frequency Born series):
\begin{equation}
\label{bornn}A(\t' , \t , k)\chi^{\ast}(k) = \q(k(\t' - \t))\chi^{\ast}(k) +
\suma_{j=2}^{\infty} \mathrm{q}_j (q)(\t' , \t , k)\, ,
\end{equation}
where
\[
\mathrm{q}_j (q)(\t',\t , k) =\chi^{\ast}(k) \int_{\Rn} e^{-ik\t' \,\cdot\, y}
(q R_k )^{j-1} (q(\cdot )e^{ik\t\cdot (\cdot)})(y)dy\, ,
\]
and $\chi^{\ast}$ is  a   cutoff function  near the origen (see the notations).

We deal with the backscattering inverse problem, for which one
assumes the data with the direction of the receiver opposed to the
source direction (echoes),  i.e. $A(-\t , \t , k )\,$. The
inverse problem is  then formally well determined.
 In this  case   the Neumann-Born series for the scattering  amplitude  is
\begin{equation}
\label{bornseries} A(-\t , \t , k )\chi^{\ast}(k) = \q(\xi)\chi^{\ast}(k) +
\suma_{j=2}^{\infty} \widehat{\widetilde{Q}_{j}(q) }(\xi)\, ,
\end{equation}
where $\xi = -2k\t $ and the $j$-adic term in the Neumann-Born
series $\widetilde{Q}_{j}(q)$ is given by the operator \eqref{qtilde}. We  define
the Born  approximation  for high frequency backscattering  data  as
\[
\widehat{q_{B,H} }(\xi)=A(-\t , \t , k )\chi^{\ast}(k)
\]
where $\xi=-2k\t\,$.
Notice that the series (\ref{bornseries}) is addapted to the reconstruction of singularities, since    $q_{B,H}-q_B$  and $q-\F^{-1}(\q(\cdot)\chi^{\ast}(|\cdot|/2))$ are $C^{\infty}$ functions.

We denote the remainder term in the  high frequency series as
$${{\bf R}}_l(q)=\sum_{j=l}^{\infty}{\widetilde{Q}_{j}(q) }.$$

The main part  of this  work, which is    \S 3,  is due to obtain  the  control  of  the  term $\widetilde{Q}_4(q)$ in dimension three:
\begin{teorema}
\label{theoremone}{\it Let us assume that $q$ is a compactly
supported function in $W^{\a\, ,\, 2}(\Rtres)\,$, for $0\leq
\a<3/2\,$. Then $\widetilde{Q}_{4}(q)\in W^{\b\, ,\, 2}(\Rtres) \,$, for any $\b$ such that $0\leq\b<\a+1/2\,$.}
\end{teorema}
We also prove in section \S 4:
\begin{teorema}
\label{lestimate}{\it Let $n\in\{2,3\}$, $q\in W^{\a , 2}(\Rn)$
compactly supported and $0\leq \a < n/2$. Assume that $C_0
> \max\{(2\norm{q}_{W^{\a , 2}})^4 , 1\} \,$, $l=4$ if $n=2$ and $l=5$ if $n=3$, see  (\ref{qtilde}). Then, for any $\b\in\R$ such that $\b<\a+1$  the remainder in  the high frequency Born series ${\bf R}_l$ converges to  a function in  $W^{\b ,2}(\Rn)$.}
\end{teorema}
 From Sobolev embedding theorem  we obtain,
 \begin{coro}\label{lcontinua} In the hypothesis of   Theorem  \ref{lestimate}, assume also    $\a>0$ in 2D and
$\a>\um$ in 3D.  Then ${\bf R}_l $ is a H\"older continuous
function.
\end{coro}

  Theorem 1  in  the case   $0\leq \alpha<n/2$  will follow  from the   above   theorems,  together  with the following estimates   for the quadratic and
cubic terms, see \cite{RV} and \cite{Re1}:
\begin{align}
\label{A}&\norm{\widetilde{Q}_2 (q)}_{\homob}\leq C\,
\norm{q}_{\homomenos}
\norm{q}_{\homoa}\, ,\\
\label{B}&\norm{\widetilde{Q}_3 (q)}_{\homob}\leq C \left(
\norm{q}_{L^2 }^2  + \norm{q}_{L^2 }\norm{q}_{\dot{W}^{-{1\over
2}-\epsi , 2} }  + \norm{q}_{\homomenos}\norm{q}_{\dot{W}^{-\epsi
, 2}} \right) \norm{q}_{\homoa}\, ,
\end{align}
where the dimension is $n=3$, $\b < \a + 1/2\,$, $\epsi
:=\a+\um-\b
>0$.  For dimension $n=2$   we have \eqref{A}   for $\b < \a + 1/2\,$  and
\begin{equation}
\label{C}\norm{\wt{Q}_3(q)}_{\homob}\leq C\, \left(\|q\|_{L^2
}\|q\|_{\homomenos} + \|q\|_{L^2 }^2 \right)\,\|q\|_{\dot{W}^{\a
\, , \, 2}}\, ,
\end{equation}
when $0\leq \b < \a +1$.

In   \S 5 we give the  procedure  to  extend the  above results  to the  case $\alpha \geq n/2$. The key  is Theorem \ref{teor_partes} which is  a  Leibniz'  type formula for derivatives  of  multiple scattering  terms.

The proofs are very involved  and  technical. For  this  reason, we only  include  the details of the  proof in  the  key  case  of Theorem \ref{theoremone}, see Proposition \ref{propprincipal} in  \S 3.1.
In other cases, we  just sketch the proof  and try to convince the  reader that similar arguments    work.

\section{Proof of Theorem \ref{theoremone}.}
Theorem \ref{theoremone} follows from the estimate
\begin{propo}{\it For $q, \a$ under the same hypothesis of Theorem
\ref{theoremone} it holds
\begin{align}
\label{D1}&\norm{\wt{Q}_4(q)}_{\homob} \leq C\,
\left(\norm{q}_{L^2}^3 + \norm{q}_{\homomenos}\norm{q}_{L^2}^2 +
\norm{q}_{\homomenos}\norm{q}_{\dot{W}^{-\epsi,2}}\norm{q}_{L^2}\right.\\
\label{D2}&\quad\quad\quad\quad\quad\quad\quad\quad \left.+\,
\norm{q}_{L^2}^2\norm{q}_{\dot{W}^{-\um-\epsi\, ,\,
2}}\right)\norm{q}_{\homoa} ,
\end{align}
for all $\b\in\R$ such that $0\leq\b < \a + 1/2\,$, where $\epsi
:=\a+\um-\b
>0$.}
\end{propo}

The quartic term in the Neumann-Born series for backscattering
data is given by
\[
\widehat{Q_{4}(q)}(\xi):=\int_{\Rtres}e^{ik\t\cdot y}(qR_{+}(k^2
))^3 (q(\cdot)e^{ik\t\cdot (\cdot )})(y)dy\, ,
\]
for any $\xi\in\Rtres\,$, where $\xi=-2k\t\,$, that is,
$k={\xim\over 2}$ and $\t=-{\xi\over \xim}\,$. From Lemma 3.1 in
\cite{R2}, this term can written as
\begin{propo}
\label{prop3.2}{\it For any dimension $n$ and
$\eta\in\Rn\setminus\{0\}\,$,
\begin{align}
\label{tpuro4}\widehat{Q_{4}(q)}(\eta)&= p.v.
\int_{\Rn}\int_{\Rn}\int_{\Rn}{\q(\xi)\q(\eta-\tau)\q(\tau-\phi)\q(\phi-\xi)\over
[\esca]\, [\tau\cdot (\eta-\tau)]\, [\phi\cdot (\eta-\phi)]}\, d\xi d\tau d\phi\\
\label{tmixto4a1}&+ 2{i\pi\over \e}\, p.v.
\int_{\Rn}\int_{\Rn}\int_{\Gae}{\q(\xi)\q(\eta-\tau)\q(\tau-\phi)\q(\phi-\xi)\over
[\tau\cdot (\eta-\tau)]\, [\phi\cdot (\eta-\phi)]}\, d\s(\xi)d\tau d\phi\\
\label{tmixto4a2}&+ {i\pi\over \e}\, p.v.
\int_{\Gae}\int_{\Rn}\int_{\Rn}{\q(\xi)\q(\eta-\tau)\q(\tau-\phi)\q(\phi-\xi)\over
[\xi\cdot (\eta-\xi)]\, [\tau\cdot (\eta-\tau)]}\, d\xi d\tau d\s(\phi)\\
\label{tmixto4b1}&- 2{\pi^2\over \e^2}\, p.v.
\int_{\Gae}\int_{\Gae}\int_{\Rn}{\q(\xi)\q(\eta-\tau)\q(\tau-\phi)\q(\phi-\xi)\over
\xi\cdot (\eta-\xi)}\, d\xi d\s(\tau) d\s(\phi)\\
\label{tmixto4b2}&-{\pi^2\over \e^2}\, p.v.
\int_{\Rn}\int_{\Gae}\int_{\Gae}{\q(\xi)\q(\eta-\tau)\q(\tau-\phi)\q(\phi-\xi)\over
\phi\cdot (\eta-\phi)}\, d\s(\xi)d\s(\tau) d\phi\\
\label{tesferico4}&- {i\pi^3 \over \e^3
}\int_{\Gae}\int_{\Gae}\int_{\Gae}\q(\xi)\q(\eta-\tau)\q(\tau-\phi)\q(\phi-\xi)\,
d\s(\xi)d\s(\tau)d\s(\phi)\, .
\end{align}}
\end{propo}

The key to understand the structure of the quartic
term is the pure spherical measures part  \eqref{tesferico4}. Hence we define,  for any
$\eta\in\Rtres\setminus\{0\}\,$

\noindent \textbf{Notation:}
 \begin{align}
\label{exp2}&\widehat{Q(q)}(\eta):={1\over \e^3
}\int_{\Gae}\int_{\Gae}\int_{\Gae}\q(\xi)\q(\eta-\tau)\q(\tau-\phi)\q(\phi-\xi)\,
d\s(\xi)d\s(\tau)d\s(\phi)\, .
\end{align}

We prove in section \S 3.1:
\begin{propo}\label{propprincipal}
{\it Let $q\in W^{\a , 2}(\Rtres)$ be a compactly supported
function with $0\leq\a < 3/2$. Then for all $\b < \a + 1/2$,
\begin{align*}
\norm{ Q(q)}_{\homob} &\leq C\, \left(\norm{q}_{L^2}^3 +
\norm{q}_{\homomenos}\norm{q}_{L^2}^2 +
\norm{q}_{\homomenos}\norm{q}_{\dot{W}^{-\epsi,2}}\norm{q}_{L^2}\right.\\
&\quad\quad\,\,\, \left.+\,
\norm{q}_{L^2}^2\norm{q}_{\dot{W}^{-\um-\epsi\, ,\,
2}}\right)\norm{q}_{\homoa}\, ,
\end{align*}
where $\epsi:= \a+\um-\b>0$ and the constant $C>0$ just depends of
$\a$, $\b$, $\epsilon$ and the support of $q$.
 }
\end{propo}
 Now  we sketch the estimates of    principal value terms
\eqref{tpuro4}-\eqref{tmixto4b2} We use a decomposition  of  the
space into diadic shelves,    as it  was done for the  cubic term
in 2D, see \cite{Re1}, and  for the  quadratic  and cubic terms in
3D in \cite{RV}. More detail can be seen in  \cite{Re2}.

Let us  state, as  a model,  the main features to control  the
principal value term $Q'(q)$, given by   \eqref{tmixto4b2},
\begin{equation}
\label{(4)}\F (Q'(q))(\eta):={1\over \e^2}\, p.v.
\int_{\Rtres}\int_{\Gae}\int_{\Gae}{\q(\xi)\q(\eta-\tau)\q(\tau-\phi)\q(\phi-\xi)\over
\phi\cdot (\eta-\phi)} \, d\s(\xi) d\s(\tau) d\phi .
\end{equation}

 The key to estimate this principal value operator is to
control the   term:
\begin{equation}
\label{not}\F ( Q_{\d}(q))(\eta):=\chi_{(\d^{-1} ,
\infty)}(\e){1\over \e^4}\,
\int_{\Gad}\int_{\Gae}\int_{\Gae}\q(\xi)\q(\eta-\tau)\q(\tau-\phi)\q(\phi-\xi)\,
d\s(\xi) d\s(\tau) d\phi ,
\end{equation}
where \begin{equation}
\label{revision}\Gad :=\{\phi\in\Rtres : ||\phi-\eta / 2|-\e / 2|\leq\d\e
\}
\end{equation}.

 Comparing  \eqref{not} with \eqref{exp2}, we observe that we
replace the sphere $\Gae$, in which the variable $\phi$  runs, by its tubular  neighborhood
  $\Gad$ of width $\d\e$ in the normal direction. Notice that $d\s_{\eta}(\phi) = \displaystyle\lim_{\d \ra 0}
{1\over \d\e}\, \chi_{\Gad}(\phi)\, d\phi$, where
$d\s_{\eta}(\phi)$ denotes the measure on the sphere $\Gae$
induced by Lebesgue measure $d\phi$.  For $\d$ small, we have that
$\wh{Q_{\d}(q)}(\eta)\sim \chi_{(\d^{-1} , \infty)}(\e)\, \d
\,\wh{Q(q)}(\eta)$. In this way we may expect
  estimates for the Sobolev norm of $Q_{\d}(q)$  obtained from estimates   of $Q(q)$ multiplied by $\d$.
If one follows  the lines of the proof of Proposition
\ref{propprincipal}, one  gets the following
\begin{lema}\label{lema}{\it
Assume that $q\in W^{\a , 2}(\Rtres)$ is compactly supported and
$0\leq \a <3/2$. Let $\b<\a+1/2$ and $\epsi=\a+1/2-\b>0$. Then
there exist $\d_1>0$ and $\g=\g(\epsilon)>1 $ such that for all
$\d$ satisfying $0<\d\leq \d_1$
 it holds
\begin{align}
\label{(2)}\norm{Q_{\d}(q)}_{\homob} &\leq C \d^{\g}\,
\left(\norm{q}_{L^2}^3 + \norm{q}_{\homomenos}\norm{q}_{L^2}^2 +
\norm{q}_{\homomenos}\norm{q}_{\dot{W}^{-\epsi,2}}\norm{q}_{L^2}\right.\\
\label{(3)}&\quad\quad\quad\,\,\, \left.+\,
\norm{q}_{L^2}^2\norm{q}_{\dot{W}^{-\um-\epsi\, ,\,
2}}\right)\norm{q}_{\homoa} ,
\end{align}
 where  $C$ just depends on $\a$, $\b$,
$\d_1$ and the support of $q$.}
\end{lema}

Now, to estimate the term \eqref{(4)}, we use a decomposition of
the Euclidean space $\Rtres $ in a similar way as was done in
\cite{Re1} for 2D:
\begin{equation}
\label{f11}\Rtres = \Gamma_{j_1^-}(\eta)\cup\bigcup_{j\,
=j_1}^{[\log_2\e ]}\Gamma_j(\eta)\cup\Gamma_{\infty}^{\ast}(\eta)
,
\end{equation}
where $j_1$ is the lowest integer such that $j_1\geq 1-\log_2 (
\d_1 )$ with $\d_1$ from Lemma \ref{lema}, $\e\geq 2^{j_1 -1}$ and
\begin{align*}
&\Gamma_{j_1^-}(\eta):=\{\phi\in\Rtres : ||\phi- \eta /
2|- \e /  2 |>2^{-j_1+1}\e\} ,\\
&\Gamma_j(\eta):=\{\phi\in\Rtres : ||\phi- \eta /  2 |- \e /
2 |\sim 2^{-j}\e\} ,\quad j_1\leq j\leq [\log_2\e] ,\\
&\Gamma_{\infty}^{\ast}(\eta):=\{\phi\in\Rtres : ||\phi- \eta / 2
|- \e /
 2|\leq 2^{-[\log_2\e ]-1}\e\} .
\end{align*}
\begin{nota}Technically this partition only makes sense for $j_1\geq
3$, but this is not a constraint   if we demand $\d_1\leq
1/4$, since $j_1\geq 1-\log_2 ( \d_1 )$. Notice that
$\Gamma_{\infty}^{\ast}(\eta)\subset \Gainf$, where
\[
\Gainf :=\{\phi\in\Rtres : ||\phi-\eta /2|-\e / 2|<1\}.
\]
\end{nota}

This decomposition  is used to  split  the operator  \eqref{(4)}.
To  control the operator  corresponding to  the annulus terms
Lemma \ref{lema}, with $\delta=2^{-j+1}$, suffices.  To deal with
  the central term, corresponding to $\Gamma_{\infty}^{\ast}(\eta)$, which is close
to the singularity $\Gae$, we use again  Lemma \ref{lema} and the
following

\begin{lema}\label{lema2}
{\it Let $f_j\in W^{\a , 2}(\Rtres)$, $j=1,\dots ,4$, be functions
such that $f_1$, $f_2$ are compactly supported and $0\leq \a<3/2$.
We denote
\begin{align}
\label{exp20} \F(Q_{\infty}^{\ast }(f_1 , f_2 , f_3 ,
f_4))(\eta)&:= \chi^{\ast}(\eta){1\over\e^3}\,
\int_{\Gainf}\int_{\Gae}\int_{\Gae}
|\wh{f}_1(\xi)\wh{f}_2(\eta-\tau)|\\
\label{exp21}&\quad\times
|\wh{f}_3(\tau-\phi)\wh{f}_4(\phi-\xi)|\, d\s(\xi) d\s(\tau) d\phi
.
\end{align}
Then for any $\b < \a + 1 /2$,
\begin{align}
\notag&\norm{Q_{\infty}^{\ast}(f_1 , f_2 , f_3 , f_4)}_{\homob}\\
\label{exp8}&\leq C(\a , \b , \supp f_1 , \supp f_2 )\, (
 \suma_{\s\in S_4} \|f_{\s (1)}\|_{L^2 }\|f_{\s (2)}\|_{L^2}\|f_{\s
 (3)}\|_{L^2
 }\|f_{\s (4)}\|_{\homoa
 }\\
\label{exp9}&\quad\,\,\quad\quad\,\,\,\quad\quad\quad\quad\quad +
\suma_{\tau\in S_4} \|f_{\tau (1)}\|_{\homomenos }\|f_{\tau
(2)}\|_{L^2}\|f_{\tau (3)}\|_{L^2 }\|f_{\tau (4)}\|_{\homoa
 }\\
\label{exp10}&\quad\,\,\quad\quad\,\,\,\quad\quad\quad\quad\quad +
 \suma_{\o\in S_4} \|f_{\o (1)}\|_{\homomenos }\|f_{\o (2)}\|_{\dot{W}^{-\epsi , 2}}\|f_{\o
 (3)}\|_{L^2 }\|f_{\o (4)}\|_{\homoa
 }\\
 \label{exp11}&\quad\,\,\quad\quad\,\,\,\quad\quad\quad\quad\quad +
 \suma_{\rho\in S_4} \|f_{\rho (1)}\|_{\dot{W}^{-\um-\epsi , 2} }\|f_{\rho (2)}\|_{L^2}\|f_{\rho
 (3)}\|_{L^2
 }\|f_{\rho (4)}\|_{\homoa
 }\, )\, ,
\end{align}
where $\epsi:=\a+\um-\b>0$.

}
\end{lema}

To estimate this central term, dealing with  the principal value, one needs to use the cancelation.
  We must replace the integral on the ring
$\Gamma_{\infty}^{\ast}(\eta)$ by
$\int_{\Gamma_{\epsi}^+(\eta)}+\int_{\Gamma_{\epsi}^-(\eta)}$,
where
\begin{align*}
&\Gamma_{\epsi}^{\,+}(\eta)=\{\xi\in \Rtres : \epsi < |\xi-\eta /2|-\e /2< 2^{-[\log_2\e]-1}\e\}\\
&\Gamma_{\epsi}^{\,-}(\eta)=\{\xi\in \Rtres : \epsi <
\e/2-|\xi-\eta / 2 |< 2^{-[\log_2\e]-1}\e \}.
\end{align*}
The map $F: \Gamma_{\epsi}^+(\eta) \to \Gamma_{\epsi}^-(\eta) $,
given by symmetry with respect to $\Gamma(\eta)$ allows us to pass
to the limit when $\epsi\ra 0^+$. To  cancel the singularities  we
use an estimate, due to Calder\'on,  for first differences in
terms of the Hardy-Littlewood maximal operator $M$ (as in \cite{Re1}, several   standard reductions are also needed):
\begin{lema}[see \cite{H}]
{\it Let $u\in W^{1, p}(\Rn), \, p>1, \, a \in \Rn$. Then
\begin{equation}
 \label{f12}
 |u(x)-u(x-a)|\leq C\, |a|\, [M(\nabla u )(x) + M(\nabla u) (x-a)].
 \end{equation}}
 \end{lema}
  After some
 changes of variables  in the integrals involving $F$,
we reduce to study the  following terms:
\begin{align*}
&{1\over\e^3}\,
\int_{\Gainf}\int_{\Gae}\int_{\Gae}|\q(\xi)\q(\eta-\tau)\q(\tau-\phi)|M\nabla\q
(\phi-\xi)\, d\s(\xi) d\s(\tau) d\phi ,\\
&{1\over\e^3}\,
\int_{\Gainf}\int_{\Gae}\int_{\Gae}|\q(\xi)\q(\eta-\tau)|M
\q(\tau-\phi)M\nabla\q
(\phi-\xi)\, d\s(\xi) d\s(\tau) d\phi ,\\
&{1\over\e^3}\,
\int_{\Gainf}\int_{\Gae}\int_{\Gae}|\q(\xi)\q(\eta-\tau)|M\nabla\q(\tau-\phi)M\q
(\phi-\xi)\, d\s(\xi) d\s(\tau) d\phi ,\\
&{1\over\e^3}\,
\int_{\Gainf}\int_{\Gae}\int_{\Gae}|\q(\xi)\q(\eta-\tau)|M\nabla
\q(\tau-\phi)|\q (\phi-\xi)|\, d\s(\xi) d\s(\tau) d\phi .
\end{align*}

The proof of Lemma \ref{lema2} follows the lines of the proof of
Proposition \ref{propprincipal}. Heuristically, Lemma \ref{lema2}
is derived from Proposition \ref{propprincipal} replacing the
domain $\Gae$ for the variable $\phi$ by the tubular neighborhood
$\Gainf$ which is the result of widening the sphere $\Gae$ a
distance  1 in the normal direction. Nevertheless there is an
additional difficulty which has to be managed: the fact that
neither $f_3$ nor $f_4$ are compactly supported and their Fourier
transform can not be controlled by the maximal operator using
Lemma \ref{remark}. But we must keep in mind that we can apply
Lemma \ref{remark} to two functions, $f_1$, $f_2$, which are
compactly supported, and the integral of $|\wh{f}_3|^2$ or
$|\wh{f}_4|^2$ in $\phi$ can be bounded by the $L^2$-norm using
that the variable $\phi$ is solid. After these comments, we omit
the long and tedious proof of Lemma \ref{lema2}. The reader can
see all the details in a similar situation for the cubic term in
2D (Lemma 2.2.3 of \cite{Re2}).

%%%%%%%%%%%%%%%%%%%%%%%%%%%%%%%%%%%%%%%%%%%%%%%%%%%%%%%%%%%%%%%%%%%%%%%%%%%%%%%%%%%%%%%%%%%
%%%%%%%%%%%%%%%%%%%%%%%%%%%%%%%%%%%%%%%%%%%%%%%%%%%%%%%%%%%%%%%%%%%%%%%%%%%%%%%%%%%%%%%%%

The key to control the principal value term \eqref{tpuro4} remains
in the following lemma whose proof follows the lines of the proof
of Proposition \ref{propprincipal}.
\begin{lema}\label{lemad}{\it We denote
\begin{align*}
\F (Q_{\d_1, \d_2, \d_3}(q))(\eta) &:= 1/\e^6\,\,
\chi_{(\d_1^{-1},\infty)}(\e)\,
\chi_{(\d_2^{-1},\infty)}(\e)\, \chi_{(\d_3^{-1},\infty)}(\e)\\
&\quad\times
\int_{\Gamma_{\d_1}(\eta)}\int_{\Gamma_{\d_2}(\eta)}\int_{\Gamma_{\d_3}(\eta)}|\q(\xi)\q(\eta-\tau)
\q(\tau-\phi)\q(\phi-\xi)|\, d\xi d\tau d\phi .
\end{align*}
Let $q$ and $\a$ as in Lemma \ref{lema} and $\b<\a+1/2$,
$\epsi=\a+1/2-\b>0$. Then there exist $\d_0>0$ and
$\g^{\ast}=\g^{\ast}(\epsilon)>1 $ such that for all $\d_1, \d_2
,\d_3$ satisfying $0<\d_1 , \d_2 , \d_3\leq \d_0$
 it holds
\begin{align*}
\norm{Q_{\d_1 , \d_2 , \d_3}(q)}_{\homob} &\leq C (\d_1 \d_2 \d_3
)^{\g^{\ast}}\, \norm{q}^4_{W^{\a , 2}(\Rtres)} ,
\end{align*}
 where  $C$ only depends on $\a$, $\b$,
$\d_0$ and the support of $q$.}
\end{lema}
Analogously to the comment about Lemma \ref{lema} above, this
result should not be surprising since $\F Q_{\d_1 , \d_2 ,
\d_3}(q)(\eta)\sim \d_1\d_2\d_3 \F Q(q)(\eta)$, for $\d_1, \d_2,
\d_3$ small. To estimate the term \eqref{tpuro4} we have to take
the partition \eqref{f11} of $\Rtres$ with $j_1$ the lowest
integer such that $j_1\geq 1-\log_2 ( \d_0 )$, for $\d_0$ from
Lemma \ref{lemad}. In particular, the control of the ring terms
\[
\int_{\Gamma_j(\eta)}\int_{\Gamma_k(\eta)}\int_{\Gamma_l(\eta)}{\q(\xi)\q(\eta-\tau)\q(\tau-\phi)
\q(\phi-\xi)\over [\xi\cdot (\eta-\xi)] [\tau\cdot(\eta-\tau)]
[\phi\cdot (\eta-\phi)]}\, d\xi d\tau d\phi \, ,\quad\quad j_1\leq
j,k,l\leq [\log_2\e]
\]
follows from Lemma \ref{lemad} with $\d_1=2^{-j+1}, \d_2 =
2^{-k+1}, \d_3 = 2^{-l+1}$, together with the fact that
\[|\xi\cdot(\eta-\xi)|= (|\xi-\frac{\eta}{2}|+|\frac{\eta}{2}|)(|\xi-\frac{\eta}{2}|-|\frac{\eta}{2}|)\geq c|\eta|^2 \delta_1,
\] where    $c>0$ and we   use  definition (\ref{revision}).

To estimate the central term
\[
\int_{\Gamma_{\infty}^{\ast}(\eta)}\int_{\Gamma_{\infty}^{\ast}(\eta)}\int_{\Gamma_{\infty}^{\ast}(\eta)}
{\q(\xi)\q(\eta-\tau)\q(\tau-\phi)
\q(\phi-\xi)\over [\xi\cdot (\eta-\xi)] [\tau\cdot(\eta-\tau)]
[\phi\cdot (\eta-\phi)]}\, d\xi d\tau d\phi
\]
we must replace each integral on the ring
$\Gamma_{\infty}^{\ast}(\eta)$ by
$\int_{\Gamma_{\epsi}^+(\eta)}+\int_{\Gamma_{\epsi}^-(\eta)}$.
 We, then, use the map $F: \Gamma_{\epsi}^+(\eta) \to \Gamma_{\epsi}^-(\eta) $.  We need again   Calder\'on estimate for  first differences  and  its  analogous estimate for second differences:
 \begin{lema}
{\it Let $u\in W^{2, p}(\Rn), \, p>1, \, a,b,c\in \Rn$. Then
\begin{equation}
 \notag |u(x-a)+u(x+b)-u(x)-u(x+b-a)|\leq C\,
|a|\,|b|\suma_{j=1}^4
M^2|D^2 u|(x_j),
\end{equation}
 where $D^2 u$ denotes the matrix of derivatives of order
 two   and $x_1=x $,
$x_2=x-a $, $x_3 = x+b$, $x_4 =x+b-a $. }
\end{lema}
These tools allow us to reduce to a sum of integrals, analogous
to those written  after Lemma 3.3  for  the case \eqref{tmixto4b2}.

 \subsection{Proof of Proposition \ref{propprincipal}.}$ $

Let us split the set $\dominio$ into the following regions
\begin{align*}
&I(\eta):=\left\{(\xi , \tau , \phi )\in\dominio \, :
\, |\phi-\xi|\geq\metac  \, ,\, |\phi-\tau|\geq {\e\over 100}  \right\}\, ,\\
&II(\eta):=\left\{(\xi , \tau , \phi )\in\dominio \, :
\, |\phi-\xi|\geq\metac\, ,\, |\phi-\tau|\leq {\e\over 100} \right\}\, ,\\
&III (\eta):=\left\{(\xi, \tau ,\phi )\in \dominio \, : \,
|\phi-\xi|\leq {\e\over 100}\, ,\, |\phi-\tau|\geq {\e\over
100}\right\}\, ,\\
&IV (\eta) :=\left\{(\xi, \tau , \phi)\in\dominio \, : \,
|\phi-\xi|\leq {\e\over 100}\, ,\, |\phi-\tau|\leq {\e\over
100}\right\}\, .
\end{align*}
In this way, we can write $Q(q)=Q_{I}(q)+ Q_{II}(q)+ Q_{III}(q)+
Q_{IV}(q)$. We will prove that
\begin{align}
\label{afir2_3d}&\norm{Q_{I}(q)}_{\homob}\leq
C\,\left(\norm{q}_{L^2}^3 +
\norm{q}_{\homomenos}\norm{q}_{L^2}^2\right)\norm{q}_{\dot{W}^{\b-1/2\,
,\, 2}}\, ,\\
\label{afir1_3d}&\norm{Q_{II}(q)}_{\homob}\leq
C\,\left(\norm{q}_{L^2}^3 +
\norm{q}_{\homomenos}\norm{q}_{L^2}^2\right)\norm{q}_{\dot{W}^{\b-1/2\,
,\, 2}}\, ,\\
\label{afir3_3d}&\norm{Q_{IV}(q)}_{\homob}\leq C\,
\left(\norm{q}_{\homomenos}\norm{q}_{\dot{W}^{-\epsi,2}}\norm{q}_{L^2}
+ \norm{q}_{L^2}^2\norm{q}_{\dot{W}^{-\um-\epsi\, ,\,
2}}\right)\norm{q}_{\homobmem}\, ,
\end{align}
provided that $\epsi >0\,$. Note that $Q_{III}(q)$ satisfies the
estimate \eqref{afir1_3d} since $Q_{II}(q)=Q_{III}(q)$.

\textbf{Proof of estimate \eqref{afir1_3d}. } Taking the change of
variable $\phi = \eta-\phi'\,$, we have
\begin{align*}
\widehat{Q_{II }(q)}(\eta) & ={1\over\e^3 }\int\int\int_{II
(\eta)}|\fourqq|\, \measuree\\
&={1\over\e^3 }\int\int\int_{\{ (\xi , \tau , \phi ): (\xi , \tau,
\eta-\phi )\in II (\eta) \}} |\fourq|\, \measure\, .
\end{align*}
We decompose
\[
\{ (\xi , \tau , \phi )\in\dominio\,:\, (\xi , \tau, \eta-\phi
)\in II (\eta) \}=\displaystyle\bigcup_{k=1}^{\infty} ( II _k
(\eta)\cup \widetilde{II}_k (\eta) ),
\]
where for any $k\in\mathbb N$, we denote
\begin{align*}
&II_k (\eta):=\left\{\etauphi\leq \metac\, ,\, \exiphi\geq
\metac\, ,\, |\phi-\xi|\sim 2^{-k}\e \, ,\, |\phi|\leq\xim
\right\} ,\\
&\widetilde{II}_k (\eta):=\left\{
 \etauphi\leq \metac\, ,\, \exiphi\geq \metac\, ,\, |\phi-\xi|\sim 2^{-k}\e \, , \, \xim\leq |\phi| \right\} ,
\end{align*}
with $(\xi, \tau , \phi )\in \Gae^3\,$. We have
\[
 \widehat{Q_{II }(q)}(\eta)  \leq \suma_{k=1}^{+\infty}
\left(\wh{ Q_{II_k}(q)}(\eta) + \wh{
Q_{\wt{II}_k}(q)}(\eta)\right) ,
\]
and then to prove \eqref{afir1_3d} we use
\[
\norm{ Q_{II}(q)}_{\homob}\leq \suma_{k=1}^{+\infty}
\left(\norm{Q_{II_k}(q)}_{\homob} +
\norm{Q_{\wt{II}_k}(q)}_{\homob}\right) .
\]

\addtocounter{equation}{-2}

For each $k\geq 1$ we claim
\begin{subequations}\begin{align}
\label{3dexp1}&\norm{Q_{II_k}(q)}_{\homob}\leq C\, 2^{-k /
2}\norm{q}_{L^2 }^3
\norm{q}_{\dot{W}^{\b-\um\, ,\, 2}}\, ,\\
\label{3dexp2}&\norm{Q_{\wt{II}_k}(q)}_{\homob}\leq C\,2^{-k/2}
\left( \norm{q}_{L^2} + \norm{q}_{\homomenos} \right)\norm{q}_{L^2
}^2 \norm{q}_{\dot{W}^{\b-\um\, ,\, 2}}\, .
\end{align}\end{subequations}

\addtocounter{equation}{1}

In the following, we use the notation in Lemma \ref{lemmasum},
which is the key of the proof of the above claims.

\textbf{Proof of claim \eqref{3dexp1}. } By Cauchy-Schwartz
inequality,
\begin{align}
\label{3dexp10}\widehat{Q_{II_k}(q)}(\eta) &\leq {1\over\e^3
}\left( \int\int\int_{II_k (\eta)} |
\q(\xi)\q(\tau+\phi-\eta)|^2 \measure\right)^{{1\over 2}}\\
\label{3dexp11}&\quad\times\left(\int\int\int_{II_k (\eta)
}|\q(\eta-\tau')\q(\eta-\phi'-\xi')|^2 \measureee \right)^{{1\over
2}} .
\end{align}
If we widen the sphere $\Gamma (\eta)$ until
$\Gamma_1(\eta):=\left\{ x\in\Rtres :
\left|\,|x-\etam|-\metam\,\right|<1 \right\}\,$, by part $(1)$ of
Lemma \ref{remark} we have
\begin{align*}
\int\int_{\Gae\times\Gae}|\q(\tau + \phi -\eta)|^2 d\s(\tau)
d\s(\phi) &\leq
C\,\int_{\Gae}\int_{\Gamma_1(\eta)}M\q(x+\phi-\eta)^2 dx \,
d\s(\phi)\\
&\leq C\,\s(\Gae)\norm{ M\q}_{L^2 }^2 \leq C\, \e^2 \norm{q}_{L^2
}^2\, ,
\end{align*}
where the last inequality follows from the boundedness of
Hardy-Littlewood maximal operator in $L^2(\Rtres)$ and Plancherel
identity, since the measure of $\Gae$ is $\pi\e^2$. In the
same way,
\begin{equation}
\label{(1)}\int_{\Gae}|\q(\eta-\tau')|^2 d\s(\tau')\leq C\,
\norm{q}_{L^2 }^2\, .
\end{equation}
If $(\xi , \tau , \phi)\in II_k(\eta)$ then $\xim\geq
2^{-k-2}\e\,$, and changing the order of integration in $\xi$ and
$\eta$ by Lemma \ref{lema2.2}, it holds
\begin{align}
\notag\norm{Q_{II_k}(q)}_{\homob}^2 &\leq C\,\norm{q}_{L^2}^4
\int_{\Rtres} \e^{2\b-4} \int_{\{\xi\in\Gae : \xim\geq
2^{-k-2}\e\}} |\q(\xi)|^2\\
\notag&\quad\times \int\int_{A_k(\eta)}
|\q(\eta-\phi'-\xi')|^2 d\s(\xi') d\s(\phi') d\s(\xi) d\eta\\
\notag& = C\, \norm{q}_{L^2 }^4 \int_{\Rtres} |\q(\xi)|^2
\int_{\{\eta\in\Lax\, :\, \xim\geq 2^{-k-2}\e\}}{\e\over\xim}\,
\e^{2\b-4}\\
\notag&\quad\times \int\int_{A_k(\eta)} |\q(\eta-\phi'-\xi')|^2
d\s(\xi')
d\s(\phi') d\s(\eta) d\xi\\
\label{(5)}&\leq C\, 2^k \norm{q}_{L^2 }^4
\int_{\Rtres}|\q(\xi)|^2 F_k(\xi)\, d\xi\leq C \, 2^{-k}
\norm{q}_{L^2}^6 \norm{q}_{\dot{W}^{\b-\um , 2}}^2\, ,
\end{align}
where the last inequality follows from part $(i)$ of Lemma
\ref{lemmasum} and $F_k(\xi)$ is defined in \eqref{3dexp52}. Also,
\begin{align}
\label{3dexp9}&A_{k}(\eta):=\left\{(\xi' ,
\phi')\in\Gae\times\Gae\, :\, |\xi'-\phi'|\leq 2^{-k+1}\e\, ,\,
|\eta-\phi'-\xi'|\geq \metac \right\}\, .
\end{align}

\qed

\textbf{Proof of claim \eqref{3dexp2}. } We take
$\wt{II}_k(\eta)=\wt{II}^1_k(\eta)\cup \wt{II}^2_k(\eta)\,$, where
\begin{align*}
&\wt{II}^1_k(\eta):=\left\{ \tupla\in \wt{II}_k(\eta) \, :\,
|\eta-\phi-\tau|\leq 2^{-k-3}\e \right\}\, ,\\
&\wt{II}^2_k(\eta):=\left\{ \tupla\in\wt{II}_k(\eta)\, :\,
|\eta-\phi-\tau|\geq 2^{-k-3}\e\right\}\, . \end{align*}

Let us start with the domain
\[
\wt{II}^1_k(\eta)=\{ (\xi , \tau ,\phi)\in\dominio :
|\eta-\phi-\tau|\leq 2^{-k-3}\e , |\eta-\xi-\tau|\geq \e / 100 ,
|\phi-\xi|\sim 2^{-k}\e , \xim\leq\phim \}.
\]
On this region $\etau\geq 2^{-k-3}\e\,$ holds. By Cauchy-Schwartz
inequality,
\begin{align}
\label{3dexp40}\widehat{Q_{\wt{II}^1_k}(q)}(\eta)&\leq {1\over\e^3
}\left( \int\int\int_{\wt{II}^1_k (\eta)} |
\q(\xi)\q(\eta-\tau)|^2 \measure\right)^{{1\over 2}}\\
\label{3dexp41}&\quad\times\left(\int\int\int_{\wt{II}^1_k (\eta)
}|\q(\tau'+\phi'-\eta)\q(\eta-\phi'-\xi')|^2 \measureee
\right)^{{1\over 2}} .
\end{align}
By part $(1)$ of Lemma \ref{remark}, we have
$\int_{\Gae}|\q(\xi)|^2 d\s(\xi)\leq C\norm{q}_{L^2}^2\,$, and by
this lemma and Fubini's theorem, for each $\xi' , \phi'\in\Gae\,$
in \eqref{3dexp41},
\[
\int_{\Gae}|\q(\tau'+\phi'-\eta)|^2 d\s(\tau')\leq
C\norm{q}_{L^2}^2\, .
\]
Taking the change $\zeta = \eta-\tau\,$, and changing the order of
integration in $\zeta$ and $\eta\,$ by Lemma \ref{lema2.2}, we may
write
\begin{align*}
\norm{Q_{\widetilde{II}^1_{k}}(q)}_{\homob}^2 &\leq C\norm{q}_{L^2
}^4 \int_{\Rtres} \e^{2\b-4}\int_{\{\zeta\in\Gae \, :\,
|\zeta|\geq 2^{-k-3}\e\}}|\q(\zeta)|^2
d\s(\zeta)\\
&\quad\times \int\int_{A_k(\eta)}|\q(\eta-\phi'-\xi')|^2
d\s(\xi')d\s(\phi')d\eta\\
&=C\norm{q}_{L^2 }^4 \int_{\Rtres} |\q(\zeta)|^2
\int_{\{\eta\in\Lambda (\zeta) :\, |\zeta|\geq 2^{-k-3}\e
\}}{\e\over
|\zeta|}\,\e^{2\b-4}\\
&\quad\times \int\int_{A_k(\eta)}|\q(\eta-\phi'-\xi')|^2
d\s(\xi')d\s(\phi')d\s(\eta)d\zeta\\
&\leq C \norm{q}_{L^2 }^4 \, 2^k \, \int_{\Rtres} |\q(\zeta)|^2
F_k(\zeta)\, d\zeta\leq C \, 2^{-k} \norm{q}_{L^2
}^6\norm{q}_{\dot{W}^{\b-\um\, ,\, 2}}^2\, ,
\end{align*}
where the last inequality follows from part $(i)$ of Lemma
\ref{lemmasum}, and $A_{k}(\eta)\,$, $F_k(\zeta) $ are defined in
\eqref{3dexp9}, \eqref{3dexp52}.

We go on with the region $\wt{II}^2_k(\eta)\,$. Let us split it as
follows: $\wt{II}^2_k(\eta)=\wt{II}^2_{k , a}(\eta)
\cup\wt{II}^2_{k , b}(\eta)\,$, where
\begin{align*}
&\wt{II}^2_{k,a}(\eta):=\left\{ \tupla\in \wt{II}^2_k(\eta) \, :\,
\etau\leq\phim \right\}\, ,\\
&\wt{II}^2_{k,b}(\eta):=\left\{ \tupla\in\wt{II}^2_k(\eta)\, :\,
\etau\geq\phim \right\} .
\end{align*}
On the region $\wt{II}^2_{k , a}(\eta)\,$, we know that if
$\xim\geq 2^{-k}\e$ we can follow the lines of the case
$II_k(\eta)\,$. So, splitting once more as $\wt{II}^2_{k ,
a}(\eta) = \wt{II}^2_{k , a,1}(\eta) \cup \wt{II}^2_{k ,
a,2}(\eta)\,$, where
\begin{align*}
&\wt{II}^2_{k,a,1}(\eta) :=\left\{ \tupla\in \wt{II}^2_{k,a}(\eta)
\, :\,
\xim\geq 2^{-k}\e \right\}\\
&\quad\quad\quad\quad\, =\{\tupla\in\dominio : 2^{-k-3}\e\leq
|\eta-\phi-\tau|\leq \e / 100 ,\, |\eta-\xi-\tau|\geq
\e / 100 ,\\
&\quad\quad\quad\quad\quad\quad\, |\phi-\xi|\sim 2^{-k}\e ,\, 2^{-k}\e\leq \xim \leq \phim ,\, \etau\leq\phim \} \\
 &\wt{II}^2_{k,a,2}(\eta) :=\left\{
\tupla\in\wt{II}^2_{k,a}(\eta)\, :\, \xim\leq 2^{-k}\e \right\}\\
&\quad\quad\quad\quad\, =\{\tupla\in\dominio : 2^{-k-3}\e\leq
|\eta-\phi-\tau|\leq \e / 100 ,\,  |\eta-\xi-\tau|\geq \e / 100
,\\
& \quad\quad\quad\quad\quad\quad\, |\phi-\xi|\sim 2^{-k}\e ,\,
\xim\leq 2^{-k}\e ,\, \xim \leq \phim ,\, \etau\leq\phim \} .
\end{align*}
We may write
\begin{equation}
\label{3dexp35}\norm{ Q_{\wt{II}^2_{k,a,1}}(q)}_{\homob}\leq C\,
2^{-k / 2}\norm{q}_{L^2 }^3 \norm{q}_{\dot{W}^{\b-\um\, ,\, 2}}\,
.
\end{equation}
In this way, we reduce to the case $\wt{II}^2_{k , a,2}(\eta)\,$,
where $\phim\leq  3\cdot 2^{-k}\e\,$ holds.

\begin{nota}
From the proof of claim 6 in \cite{RV} one deduces to the following
estimates:
\begin{align}
\label{f5}&\int\int_{\Gae^2 \cap \{|x|\leq |y|\leq 2^{-k+1}\e ,
|x-y|\leq \metac\}}|\q(x-y)|^2 d\s(y) d\s(x)\leq C\, 2^{-k} \e^2
\norm{q}_{\homomenos}^2 ,
\\
\label{f6}&\int\int_{\Gae^2 \cap \{|x|\leq |y|\leq 2^{-k+1}\e , \,
2^{-k-3}\e\leq |x-y|\leq \metac\}}|\q(x-y)|^2 d\s(y) d\s(x)\leq
C\, \e\, \norm{q}_{L^2 }^2 .
\end{align}
\end{nota}

By Cauchy-Schwartz inequality as in
\eqref{3dexp10},\eqref{3dexp11}, the fact that Lemma \ref{remark}
implies \eqref{(1)}, applying the change $\zeta = \eta-\tau\,$,
since
\[
\int\int_{\Gae^2 \cap \{|\zeta|\leq |\phi|\leq 3\cdot 2^{-k}\e ,
\, 2^{-k-3}\e\leq |\zeta-\phi|\leq \metac\}}|\q(\phi-\zeta)|^2
d\s(\zeta) d\s(\phi)\leq C\, \e\, \norm{q}_{L^2 }^2
\]
by \eqref{f6}, and changing the order of integration in $\xi$ and
$\eta$ by Lemma \ref{lema2.2}, we have
\begin{align}
\label{3dexp29}\norm{Q_{\wt{II}^2_{k,a,2}}(q)}_{\homob}^2 &\leq
C\norm{q}_{L^2 }^4 \int_{\Rtres} \e^{2\b-5}
\int_{\Gae}|\q(\xi)|^2 d\s(\xi)\\
\label{3dexp32}&\quad\times\int\int_{A_k(\eta)}|\q(\eta-\phi'-\xi')|^2
d\s(\xi')d\s(\phi')d\eta\\
\label{3dexp33}&\leq C \norm{q}_{L^2 }^4 \int_{\Rtres}
{|\q(\xi)|^2\over\xim }\, F_k(\xi)\, d\xi\\
\notag&\leq C \, 2^{-2k}\norm{q}_{\homomenos}^2 \norm{q}_{L^2 }^4
\norm{q}_{\dot{W}^{\b-\um\, ,\, 2}}^2\, ,
\end{align}
where the last inequality follows from part $(i)$ of Lemma
\ref{lemmasum} and $A_{k}(\eta)\,$, $F_k(\xi) $ are defined in
\eqref{3dexp9}, \eqref{3dexp52}.

Let us go on with the domain
\begin{align*}
\wt{II}^2_{k,b}(\eta) &=\{\tupla\in\dominio : 2^{-k-3}\e\leq
|\eta-\phi-\tau|\leq \e / 100 ,\\
& \quad\quad  |\eta-\xi-\tau|\geq \e / 100 ,\, |\phi-\xi|\sim
2^{-k}\e ,\, \xim\leq \phim \leq \etau \} .
\end{align*}
So, if $\tupla\in \wt{II}^2_{k , b}(\eta)\,$ then $\etau\geq
2^{-k-4}\e\,$. Following the lines for the case
$\wt{II}^1_k(\eta)\,$, one obtains
\begin{equation}
\label{3dexp37}\norm{Q_{\wt{II}^2_{k , b}}(q)}_{\homob}\leq C \,
2^{-k/2} \norm{q}_{L^2 }^3\norm{q}_{\dot{W}^{\b-\um\, ,\, 2}}\, .
\end{equation}

We conclude the estimate
\begin{equation}
\label{3dexp38}\norm{Q_{\wt{II}^2_{k}}(q)}_{\homob}\leq C \,
2^{-k/2} \left( \norm{q}_{L^2 } + \norm{q}_{\homomenos}
\right)\norm{q}_{L^2 }^2 \norm{q}_{\dot{W}^{\b - \um\, ,\, 2}}\, ,
\end{equation}
and the claim \eqref{3dexp2}.

This ends the proof of estimate \eqref{afir1_3d}.

\qed

\proofof{estimate \eqref{afir2_3d} }Taking the change of variable
$\phi = \eta-\phi'\,$, we have
\begin{align*}
\widehat{Q_{I }(q)}(\eta) & ={1\over\e^3 }\int\int\int_{I
(\eta)}|\fourqq|\, \measuree\\
&={1\over\e^3 }\int\int\int_{\{ (\xi , \tau, \eta-\phi )\in I
(\eta) \}} |\fourq|\, \measure\, .
\end{align*}

For $\eta\in\Rtres $ fixed, we take the decomposition
\[
\left\{ (\xi, \tau , \phi )\in \Gae^3\, :\, |\eta-\phi-\xi|\geq
\metac , |\eta-\phi-\tau|\geq\metac\right\}= I_1 (\eta)\cup
\widetilde{I}_1 (\eta) \cup I _{2}(\eta)\, ,
\]
where
\begin{align*}
&I_1 (\eta):=\left\{|\phi-\tau|\leq {\e\over 400}\, ,\,
\phim\leq\etau\, ,\, \etauphi\geq \metac\, ,\,  \exiphi\geq \metac
\right\}\, ,\\
&\widetilde{I}_1 (\eta):=\left\{|\phi-\tau|\leq {\e\over 400}\,
,\, \phim\geq\etau\, ,\, \etauphi\geq \metac\, ,\,  \exiphi\geq
\metac
\right\}\, ,\\
&I_{2}(\eta):=\left\{|\phi-\tau|\geq {\e\over 400}\, ,\,
\etauphi\geq \metac\, ,\,  \exiphi\geq \metac \right\}\, ,
\end{align*}
with $(\xi, \tau , \phi )\in \Gae^3\,$. It holds
\[
 \widehat{Q_{I }(q)}(\eta) = \wh{Q_{I_1}(q)}(\eta) +
\wh{Q_{\wt{I}_1}(q)}(\eta) + \wh{Q_{I_{2}}(q)}(\eta)\, ,
\]
and also,
\[
\norm{ Q_{I}(q)}_{\homob}\leq \norm{Q_{I_1}(q)}_{\homob} +
\norm{Q_{\wt{I}_1}(q)}_{\homob} + \norm{Q_{I_{2}}(q)}_{\homob}\, .
\]
We claim the following:\addtocounter{equation}{-18}
\begin{subequations}\begin{align}
\label{3dexp43}&\norm{Q_{I_1}(q)}_{\homob}\leq C\, \norm{q}_{L^2
}^3
\norm{q}_{\dot{W}^{\b-\um\, ,\, 2}}\, ,\\
\label{3dexp44}&\norm{Q_{\wt{I}_1}(q)}_{\homob}\leq C\,
\norm{q}_{L^2 }^3
\norm{q}_{\dot{W}^{\b-\um\, ,\, 2}}\, ,\\
\label{3dexp45}&\norm{Q_{I_{2}}(q)}_{\homob}\leq C\,
\norm{q}_{\homomenos}\norm{q}_{L^2}^2 \norm{q}_{\dot{W}^{\b-\um ,
2}}\,  .
\end{align}
\end{subequations}
\addtocounter{equation}{17}

The estimate \eqref{afir2_3d} follows from these three claims. In
their proofs we use the notation introduced in the key Lemma
\ref{lemmasum} located in the appendix.

\qed

\proofof{claim \eqref{3dexp43} } On this region we have $\etau\geq
{\e\over 200}\,$. Applying the Cauchy-Schwartz inequality as in
\eqref{3dexp40}-\eqref{3dexp41}, since for
$\eta\in\Rtres$, $\phi'\in\Gae$ fixed, by Lemma \ref{remark} it
holds $ \int_{\Gae}|\q(\xi)|^2d\s(\xi)\leq C\,\norm{q}_{L^2 }^2$,
$ \int_{\Gae}|\q(\eta-\phi'-\xi')|^2 d\s(\xi')\leq
C\,\norm{q}_{L^2 }^2$ and $\s(\Gae)=\pi\e^2$, we obtain
\begin{align*}
\wh{Q_{I_1 }(q)}(\eta)&\leq {C\over \e^2 }\, \norm{q}_{L^2 }^2 \,
\bigg(\int_{\{\tau\in\Gae : \etau\geq
{\e\over 200}\}}|\expdos|^2 d\s(\tau)\bigg)^{{1\over 2}}\\
&\quad\times \bigg(\int\int_{A(\eta)}|\qephitau|^2
d\s(\tau')d\s(\phi')\bigg)^{{1\over 2}}\, ,
\end{align*}
where
\begin{align}
\label{3dexp46}&A(\eta):=\bigg\{ (\tau' , \phi')\in\Gae^2 \, :\,
|\tau'-\phi'|\leq {\e\over 400}\, ,\, |\eta-\tau'-\phi'|\geq\metac
\bigg\}\, .
\end{align}
Taking $\zeta=\eta-\tau\,$ and changing the order of integration
in $\zeta$ and $\eta$ by Lemma \ref{lema2.2}, we have
\begin{align*}
\norm{Q_{I_1}(q)}_{\homob}^2 &\leq C\norm{q}_{L^2 }^4
\int_{\Rtres} |\q(\zeta)|^2 \int_{\{\eta\in\Lambda(\zeta) : \,
|\zeta|\geq
{\e\over 200}\}}{\e\over |\zeta|}\, \e^{2\b-4}\\
&\quad\times \int\int_{A(\eta)}|\qephitau|^2
d\s(\tau')d\s(\phi')d\s(\eta)d\zeta\\
&\leq C\norm{q}_{L^2 }^4 \int_{\Rtres} |\q(\zeta)|^2 F_1(\zeta)\,
d\zeta \leq C\norm{q}_{L^2 }^6 \norm{q}_{\dot{W}^{\b-\um\, ,\,
2}}^2\, ,
\end{align*}
where the last inequality follows from part $(i)$ of Lemma
\ref{lemmasum} with $k=1\,$ and $F_k(\zeta)$ is defined in
\eqref{3dexp52}. Let us remark that $A(\eta)\subset A_1(\eta)\,$
according  to the notation in \eqref{3dexp9}.

\qed

\proofof{claim \eqref{3dexp44} }We consider the
partition
$\wt{I}_1(\eta)=\wt{I}_{1,a}(\eta)\cup\wt{I}_{1,b}(\eta)\,$, where
\begin{align*}
&\wt{I}_{1,a}(\eta):=\left\{\tupla\in \wt{I}_1(\eta)\, :\,
|\eta-\phi|\leq {\e\over 200}\right\}\, ,\\
&\wt{I}_{1,b}(\eta):=\left\{\tupla\in \wt{I}_1(\eta)\, :\,
|\eta-\phi|\geq {\e\over 200}\right\}\, .
\end{align*}

On the region $\wt{I}_{1,a}(\eta)$, $\xim\geq {\e\over 200}\,$
holds. We may follow the same lines as in the proof for the case
$II_k(\eta)\,$ with $k=1$, interchanging the roles of the factors
$\q(\eta-\phi-\xi)$ and $\q(\tau+\phi-\eta)\,$. As we mentioned
before, $A(\eta)\subset A_1(\eta)\,$, according  to the notation
in \eqref{3dexp9} and we may apply part $(i)$ of Lemma
\ref{lemmasum} for $k=1\,$. We have
\[
\norm{Q_{\wt{I}_{1,a} }(q)}_{\homob}\leq C \norm{q}_{L^2 }^3
\norm{q}_{\dot{W}^{\b-{1\over 2}\, ,\, 2}}\, .
\]

On the region $\wt{I}_{1,b}(\eta)$, $\etau\geq {\e\over 400}\,$
holds. Hence, the proof of the estimate for the region $I_1(\eta)$
is valid here, deducing
\[
\norm{Q_{\wt{I}_{1,b} }}_{\homob}\leq C \norm{q}_{L^2 }^3
\norm{q}_{\dot{W}^{\b-{1\over 2}\, ,\, 2}}\, .
\]

\qed

\proofof{claim \eqref{3dexp45} }

Let $\eta\in\Rtres\setminus \{0\}\,$.  The occurrences of $q$ in
the term $\wh{Q_{I_2}(q)}(\eta)$ interact with each other, so that
we only may bound by the  maximal operator just once.
 We need to consider an extra splitting carried out by taking the  set of points $\{\t_j\, :\,
1\leq j\leq \mathcal{N}\}$ in the unitary sphere $\mathbb S^2$,
where $\mathcal{N}$ is large enough, to get a covering of the
sphere $\Gae $ with \label{constante}$\mathcal{N}$
spherical cups $J_j(\eta)\,$ centered at $\Omega_j$ of radius
$C_1\e\,$, for a certain constant $C_1>0$ to be chosen later (with
$\mathcal{N}\sim {1\over C_1^2}$ ). We define  for every  $j$
\begin{equation}
\label{tj}\Omega_j = {\eta\over 2}+\metam \,\t_j\, .
\end{equation}
Then
\begin{align*}
\Gae =\union_{j=1}^{\mathcal{N}}J_j(\eta)\, ,\quad
\wh{Q_{I_2}(q)}(\eta)\leq \suma_{j=1}^\mathcal{N}
\wh{R_{J_j}(q)}(\eta)\, ,
\end{align*}
where
\begin{align*}
\wh{R_{J_j}(q)}(\eta)&:={1\over\e^3}
\int_{J_j(\eta)}\int_{X_{\eta}(\tau)}\int_{Y_{\eta}(\phi)}
|\q(\xi)\q(\eta-\tau)\q(\tau+\phi-\eta)|\\
&\quad\times |\q(\eta-\phi-\xi)|\,d\s(\xi)d\s(\phi)d\s(\tau) ,
\end{align*}
and
\begin{align}
\label{3dexp47} &X_{\eta}(\tau):=\{\phi\in\Gae \, :\,
|\phi-\tau|\geq \e / 400\, ,\,
|\eta-\phi-\tau|\geq \e / 100 \}\, ,\\
\label{3dexp48}&Y_{\eta}(\phi):=\{\xi\in\Gae \, :\,
|\eta-\phi-\xi|\geq \e / 100 \}\, .
\end{align}

We fix $j\in\{1 , \, \dots\, , \mathcal{N}\}$. In this part, we
take an orthonormal reference of $\Rtres$ $\{e_1 , e_2 , e_3\}$
such that $e_1 = \t_j$, according to the notation used in
\eqref{tj}. On the integral expression $\wh{R_{J_j}(q)}(\eta)$, we
apply Cauchy-Schwartz inequality as in
\eqref{3dexp10}-\eqref{3dexp11}. For each $j ,\eta$ fixed, by
Fubini's theorem, we have
\begin{align*}
&\int_{J_j(\eta)}\int_{X_{\eta}(\tau)}\int_{Y_{\eta}(\phi)}|\q(\xi)\q(\tau+\phi-\eta)|^2
d\s(\xi) d\s(\phi) d\s(\tau)\\
&\leq \int_{\Gae}|\q(\xi)|^2
\int_{J_j(\eta)}\int_{X_{\eta}(\tau)}|\q(\tau+\phi-\eta)|^2
d\s(\phi) d\s(\tau) d\s(\xi)\, .
\end{align*}
Moreover,
\begin{align*}
&\int_{J_j(\eta)}\int_{X_{\eta}(\tau')}\int_{Y_{\eta}(\phi')}|\q(\eta-\tau')\q(\eta-\phi'-\xi')|^2
d\s(\xi') d\s(\phi') d\s(\tau')\\
&\leq \int_{\Gae}\int_{\Gae} \int_{Y_{\eta}(\phi')}
|\q(\eta-\tau')\q(\eta-\phi'-\xi')|^2
d\s(\xi') d\s(\phi') d\s(\tau')\\
&= \int_{\Gae} |\q(\eta-\tau')|^2 d\s(\tau')\, \int_{\Gae}
\int_{Y_{\eta}(\phi')}|\q(\eta-\phi'-\xi')|^2 d\s(\xi')
d\s(\phi')\\
&\leq C\, \norm{q}_{L^2}^2 \int_{\Gae} \int_{Y_{\eta}(\phi')}
|\q(\eta-\phi'-\xi')|^2 d\s(\xi') d\s(\phi') ,
\end{align*}
bounding by the maximal operator using Lemma \ref{remark} in the
last inequality. Changing the order of integration in $\xi$ and
$\eta$ by Lemma \ref{lema2.2}, we may write
\begin{align*}
\norm{R_{J_j}(q)}_{\homob}^2 &\leq C\, \norm{q}_{L^2 }^2
\int_{\Rtres} \e^{2\b-6} \int_{\Gae} |\q(\xi)|^2 \int_{J_j(\eta)}
\int_{X_{\eta }(\tau)} |\q(\tau+\phi-\eta )|^2 d\s(\phi)
d\s(\tau)\\
&\quad\times \int_{\Gae} \int_{Y_{\eta}(\phi')}
|\q(\eta-\phi'-\xi')|^2 d\s(\xi') d\s(\phi') d\s(\xi) d\eta\\
&=C\, \norm{q}_{L^2}^2 \int_{\Rtres} {|\q(\xi)|^2 \over\xim}\,
\int_{\Lax} \e^{2\b-5} \int_{J_j(\eta)}\int_{X_{\eta}(\tau)}
|\q(\tau + \phi -\eta)|^2 d\s(\phi) d\s(\tau)\\
&\quad\times \int_{\Gae} \int_{Y_{\eta}(\phi')}
|\q(\eta-\phi'-\xi')|^2 d\s(\xi') d\s(\phi') d\s(\eta) d\xi .
\end{align*}
We fix $\xi\in\Rtres\setminus\{0\}$ and denote
\begin{align*}
G_j(\xi) &:= \int_{\Lax} \e^{2\b-5}
\int_{J_j(\eta)}\int_{X_{\eta}(\tau)} |\q(\tau + \phi -\eta)|^2
d\s(\phi) d\s(\tau)\\
&\quad\times \int_{\Gae} \int_{Y_{\eta}(\phi')}
|\q(\eta-\phi'-\xi')|^2 d\s(\xi') d\s(\phi') d\s(\eta) .
\end{align*}
We write $\eta\in\Lax$ in cylindric coordinates $\eta = \xi + s z
$, with $s\geq 0$ and $z\in\{\xi\}^{\bot}$, $|z|=1$. It is true
that $d\s(\eta) = s \, ds \, d\s(z)$. Let $h(s):=\e = (\xim^2 +
s^2 )^{\um}$. We have
\begin{align}
\notag G_j(\xi) &= \int_{\mathbb S^1}\int_{0}^{\infty}
h(s)^{2\b-5} \int_{J_j(\xi +s z)} \int_{X_{\xi + s z}(\tau)}
|\q(\tau + \phi - (\xi + s
z))|^2 d\s(\phi) d\s(\tau)\\
\notag &\quad\times \int_{\Gamma (\xi + s z)}\int_{Y_{\xi+s
z}(\phi')} |\q(\xi + s z -\phi' -\xi ')|^2 d\s(\xi') d\s(\phi' )
\, s\, ds\,
d\s(z)\\
\label{exp13}&\leq C\, \norm{q}_{L^2}^2 \int_{\mathbb S^1}
\int_0^{\infty}
h(s)^{2\b-4}\\
\label{exp14}&\quad\times \int_{\Gamma (\xi + s z)} \int_{Y_{\xi +
sz }(\phi')} |\q(\xi + s z -\phi' -\xi ')|^2 d\s(\xi') d\s(\phi' )
\, s\, ds\, d\s(z)\, ,
\end{align}
where the last inequality follows from the following
\begin{afir}
\label{afir_nueva}{\it Let $1\leq j \leq \mathcal{N}$,
$\xi\in\Rtres$, $z\in\{\xi\}^{\bot}$, with $|z|=1$ and $s\geq 0$.
Hence
\begin{align}
\label{exp12}&\int_{J_j(\xi + sz)}\int_{X_{\xi + s z}(\tau)}
|\q(\tau + \phi - (\xi + s
z))|^2 d\s(\phi) d\s(\tau)\\
\notag&\leq C\, h(s)\norm{q}_{L^2}^2 .
\end{align} }
\end{afir}

\proofof{claim \ref{afir_nueva}}

We write $\tau$, $\phi$ in spherical coordinates with respect to
the reference $\{e_1 , e_2 , e_3\}$:
\begin{align}
\label{exp15}&\tau ={\xi + s z\over 2} + {h(s)\over 2}\, (\sin
\psi \cos\d \, e_1
+ \sin \psi \sin\d \, e_2 + \cos \psi \, e_3)\, ,\\
\label{exp16}&\phi = {\xi + s z\over 2} + {h(s)\over 2}\,
(\sin\zeta\cos\g \, e_1 + \sin\zeta\sin\g \, e_2 + \cos\zeta \,
e_3),
\end{align}
where $\psi , \zeta\in [0,\pi]$, $\d , \g\in (-\pi , \pi ]$. It
holds
\[
d\s(\phi)d\s(\tau) = h(s)^4 \sin\psi \sin\zeta\, d\g \, d\zeta \,
d\d \, d \psi .
\]
Notice that if $\tau\in J_j(\xi + s z)$ then $\tau$ belongs to the
``curvilinear square'' from the sphere $\Gae$ which contains the
spherical cup $J_j(\xi + s z)$ given by
\[
(\psi , \d)\in \left[{\pi\over 2}-\epsi_0 \, , \, {\pi\over
2}+\epsi_0 \right]\times [-\epsi_0 , \epsi_0],
\]
where $\epsi_0 = \epsi_0 (C_1)$ satisfies $\sin \epsi_0 = 2 C_1$.
For each $\zeta \in [0, \pi]$ and $\psi , \d$, we define
\begin{align*}
X^{\ast}(\zeta , \psi , \d )&:=\{\g\in (-\pi , \pi ]\, :\, \phi\in
X_{\xi + s z}(\tau)\}\\
&=\left\{\g\in (-\pi , \pi]\, :\, -(1-{1\over 5000})\leq \sin\psi
\cos\d \sin \zeta \cos\g\right.\\
&\left.\quad\,\,\,\,\, + \sin\psi \sin\d \sin\zeta \sin\g +
\cos\psi \cos \zeta \leq 1-{1\over 80000}\right\} .
\end{align*}
The integral expression \eqref{exp12} is bounded by
\begin{align*}
&C\, \int_{{\pi\over 2}-\epsi_0 }^{{\pi\over 2}+\epsi_0
}\int_{-\epsi_0 }^{\epsi_0 }\int_0^{\pi}\int_{X^{\ast}(\zeta ,
\psi , \d)}h(s)^4 \sin\psi \sin\zeta \, |\q(A(j,s,\psi , \d ,
\zeta , \g))|^2 d\g \, d\zeta \, d\d \, d\psi\, ,
\end{align*}
where
\begin{align*}
A(j,s,\psi , \d , \zeta , \g)&:={h(s)\over 2}\, ((\sin\psi \,
\cos\d \, + \, \sin\zeta \, \cos\g )\, e_1\\
&\quad + \,(\sin\psi \,\sin \d \, +\, \sin\zeta \, \sin\g )\, e_2
\, + \, ( \cos\psi \, +\, \cos\zeta ) \, e_3).
\end{align*}

For technical reasons we divide the domain which corresponds to
the angles $\psi , \d , \zeta , \g\,$ into two pieces
$\mathcal{A}_1 \, , \, \mathcal{A}_2\,$:
\begin{align*}
&\mathcal{A}_1 :=\{(\psi , \d , \zeta , \g )\, :\,
|\cos(\psi-\zeta)|<1-10^{-9} \}\, ,\\
&\mathcal{A}_2 :=\{(\psi , \d , \zeta , \g )\, :\,
|\cos(\psi-\zeta)|\geq 1-10^{-9} \}\, ,
\end{align*}
where $(\psi , \d , \zeta )\in  \left[\, {\pi\over 2}-\epsi_0 \, ,
{\pi\over 2}+\epsi_0 \right]\times [-\epsi_0 , \epsi_0 ]\times [0,
\pi ]$, and for each $\psi$, $\d$, $\zeta$ fixed $\g$ belongs to
the set $X^{\ast}(\zeta , \psi ,\d)$. In this way, \eqref{exp12}
becomes bounded by
\[
C\,\left(\int\int\int\int_{\mathcal{A}_1}+\int\int\int\int_{\mathcal{A}_2}\right)
h(s)^4 \sin\psi \sin\zeta \, |\q(A(j,s,\psi , \d , \zeta , \g))|^2
d\g \, d\zeta \, d\d \, d\psi\, .
\]

\begin{nota}
We choose $\mathcal{N}$ large enough in order to take the radius
of the spherical cup $J_j(\eta)$ with $C_1<{10^{-5}\over 2}\,$.
\end{nota}

{\it Estimate for the domain $\mathcal{A}_1$.}

\vskip 0.3cm

If $C_1<{10^{-5}\over 2}\,$ then $|\cos\psi|\leq 2C_1 <10^{-5}\,$.
We have
\begin{align*}
1-10^{-9} &> |\cos(\psi-\zeta)|\geq |\sin\psi\,
\sin\zeta\,|-|\cos\psi\, \cos\zeta\,|\\
&\geq
\sqrt{1-10^{-10}}\,\sin\zeta\, - 10^{-5}|\cos\zeta\,|\\
&= \sqrt{1-10^{-10}}\,\sqrt{1-\cos^2 \zeta }\, -
10^{-5}|\cos\zeta\,|\, ,
\end{align*}
and
\[
\cos^2 \zeta + 2\cdot 10^{-5}(1-10^{-9})|\cos\zeta\,|+
(1-10^{-9})^2 - (1-10^{-10})> 0\, ,
\]
where the discriminant of the corresponding quadratic polynomial
in $|\cos\zeta\,|$ is
\[
\Delta = 4\cdot 10^{-9} (1-10^{-10})(2-10^{-9})>0\, ,
\]
in such a way that $\zeta$ satisfies
\[
|\cos\zeta\,|>
-(1-10^{-9})10^{-5}+\sqrt{10^{-9}(1-10^{-10})(2-10^{-9})}> 3\cdot
10^{-5}\, .
\]

Keeping in mind that $|\cos\psi |<10^{-5}\, ,\,
|\cos\zeta\,|>3\cdot 10^{-5} \,$, it is true that
\begin{align}
\label{form1}&|\sin\zeta\,\cos\psi\,\cos (\d-\g) - \sin\psi\,\cos\zeta|\\
\label{form2}&\geq|\sin\psi\, \cos\zeta
|-|\sin\zeta\,\cos\psi\,\cos (\d-\g)|\geq|\sin\psi\,
\cos\zeta |-|\sin\zeta\,\cos\psi |\\
\label{form3}&\geq 3\cdot 10^{-5}\sqrt{1-10^{-10}}\, - 10^{-5}
\sqrt{1-9\cdot 10^{-10}} \sim 2\cdot 10^{-5}\, .
\end{align}

For each $j , s , \d $ fixed we take the change of variables
$(\zeta\, ,\, \gamma\, ,\, \psi)\, \ra\, \l =(\l_1 , \l_2 , \l_3
)\,$, given by
\begin{align}
\label{change1}\l &= {h(s)\over 2}\, ((\sin\psi \,
\cos\d \, + \, \sin\zeta \, \cos\g )\, e_1\\
\label{change2}&\quad + \,(\sin\psi \,\sin \d \, +\, \sin\zeta \,
\sin\g )\, e_2 \, + \, ( \cos\psi \, +\, \cos\zeta ) \, e_3).
\end{align}
We have
\[
\left| { \partial\left( \l_1\, ,\, \l_2\, ,\, \l_3  \right) \over
\partial \left( \zeta\, ,\, \g\, ,\, \psi \right) }\right| =
{h(s)^3 \sin\zeta\, \left| \sin\zeta \,\cos\psi\, \cos (\d-\gamma
) - \sin\psi\, \cos\zeta \right| \over 8}\, .
\]

By Fubini's and Toneli's theorems and estimates
\eqref{form1},\eqref{form2},\eqref{form3}, we may write
\begin{align*}
&\int\int\int\int_{\mathcal{A}_1} h(s)^4 \sin\psi \sin\zeta \,
|\q(A(j,s,\psi , \d , \zeta , \g))|^2 d\g \, d\zeta \, d\d \,
d\psi\\
&\leq C\, \int_{-\epsi_0}^{\epsi_0}\int_{\Rtres} h(s)|\q(\l)|^2
d\l \,
d\d\\
&=C\, h(s)\int_{\Rtres}|\q(\l)|^2 d\l\, .
\end{align*}

\vskip 0.3cm

{\it Estimate for the domain $\mathcal{A}_2$.}

\vskip 0.3cm

Now we apply the change of variables $(\psi\, ,\, \d \, , \,
\gamma )\, \ra\, \l =(\l_1 , \l_2 , \l_3 )\,$, given by
\eqref{change1}-\eqref{change2} for each $j,s,\zeta$ fixed. It
holds
\[
\left| { \partial\left( \l_1\, ,\, \l_2\, ,\, \l_3  \right) \over
\partial \left( \psi\, ,\, \d\, ,\, \gamma \right) }\right| =
{h(s)^3 \sin^2\psi\,\sin\zeta\, \left| \sin (\d-\gamma ) \right|
\over 8}\, .
\]
On one hand, since $|\psi-{\pi\over 2}|\leq \epsi_0$, we have that
$|\cos\psi|\leq 2C_1$. Hence the sinus of $\psi$ is lower bounded
by a strictly positive constant.

On the other hand, since $\g\in X^{\ast}(\zeta , \psi , \d )$ we
know that the expression
\begin{align*}
&\sin\psi \cos\d \sin \zeta \cos\g + \sin\psi \sin\d \sin\zeta
\sin\g + \cos\psi \cos \zeta\\
&= \cos (\psi-\zeta)\cos(\d-\g) + \cos\psi\, \cos\zeta\, (1-\cos
(\d-\g))
\end{align*}
takes values between $-(1-1/ 5000)$ and $1-1/ 80000$. We have
\begin{align*}
1-{1\over 80000} &\geq |\cos (\psi-\zeta)\cos(\d-\g) + \cos\psi\,
\cos\zeta\, (1-\cos (\d-\g))|\\
&\geq |\cos (\psi-\zeta)\cos(\d-\g)| - |\cos\psi\, \cos\zeta| \,
(1-\cos (\d-\g))\\
&\geq (1-10^{-9})|\cos(\d-\g)| - 2C_1 \, |\cos\zeta| \, (1-\cos
(\d-\g)) .
\end{align*}
Provided that $\cos(\d-\g)\geq 0$ it holds
\[
1-{1\over 80000} \geq (1-10^{-9}+2C_1 |\cos\zeta|)\cos(\d-\g) - 2
C_1 |\cos\zeta| ,
\]
and hence,
\[
\cos (\d-\g)\leq {1-1/ 80000 + 2C_1 |\cos\zeta|\over
1-10^{-9}+2C_1 |\cos\zeta|}\leq  {1-1/ 80000 + 2C_1 \over
1-10^{-9}} .
\]
Nevertheless, if $\cos(\d-\g)<0$ we may write
\begin{align*}
1-{1\over 80000} &\geq -(1-10^{-9}-2 C_1 |\cos\zeta|)\cos(\d-\g) -
2 C_1 |\cos\zeta|\\
&\geq -(1-10^{-9}-2 C_1 )\cos(\d-\g) - 2 C_1 |\cos\zeta| .
\end{align*}
 The choice $C_1 < {10^{-5}\over 2}$ allows us to  write
\[
\cos (\d-\g)\geq -{1-1/ 80000 + 2C_1 |\cos\zeta|\over 1-10^{-9} -
2 C_1 }\geq -{1-1/ 80000 + 2C_1 \over 1-10^{-9} - 2 C_1} ,
\]
and
\[
|\cos (\d-\g)|\leq {1-1/ 80000 + 2C_1 \over 1-10^{-9} - 2 C_1} .
\]
Hence   $|\sin (\d-\g)|$ is  bounded below by a strictly
positive constant which only depends on $C_1$.

By Fubini's and Toneli's theorems, we have
\begin{align*}
&\int\int\int\int_{\mathcal{A}_2} h(s)^4 \sin\psi \sin\zeta \,
|\q(A(j,s,\psi , \d , \zeta , \g))|^2 d\g \, d\zeta \, d\d \,
d\psi\\
&\leq C\, \int_{0}^{\pi}\int_{\Rtres} h(s)|\q(\l)|^2 d\l \,
d\zeta \\
&=C\, h(s)\int_{\Rtres}|\q(\l)|^2 d\l\, ,
\end{align*}
and claim \ref{afir_nueva} follows.

\qed

We return to the expression \eqref{exp13}-\eqref{exp14}. In
\eqref{exp14} we write the variables $\xi' , \phi'$ in spherical
coordinates as we did in \eqref{exp15}, \eqref{exp16}:
\begin{align*}
&\xi' ={\xi + s z\over 2} + {h(s)\over 2}\, (\sin \Theta \cos\t \,
e_1
+ \sin \Theta \sin\t \, e_2 + \cos \Theta \, e_3)\, ,\\
&\phi' = {\xi + s z\over 2} + {h(s)\over 2}\, (\sin\zeta'\cos\g'
\, e_1 + \sin\zeta'\sin\g' \, e_2 + \cos\zeta' \, e_3),
\end{align*}
where $\Theta , \zeta'\in [0,\pi]$, $\t , \g'\in (-\pi , \pi ]$.
It holds
\[
d\s(\xi')d\s(\phi') = h(s)^4 \sin\Theta \sin\zeta'\, d\t \,
d\Theta \, d\g' \, d \zeta' .
\]

For each $\Theta \in [0, \pi]$ and $\zeta' , \g' $, we define
\begin{align}
\notag Y^{\ast}(\Theta , \zeta' , \g' )&:=\{\t\in (-\pi , \pi ]\,
:\, \xi'\in
Y_{\xi + s z}(\phi')\}\\
\label{(11a)}&=\big\{\t\in (-\pi , \pi ]\, :\, \sin\Theta
\cos\t \sin \zeta' \cos\g'\\
\label{(11b)}&\quad\,\,\,\,\, + \sin\Theta \sin\t \sin\zeta'
\sin\g' + \cos\Theta \cos \zeta' \geq -\big(1-1/ 5000 \big) \big\}
.
\end{align}
\eqref{exp13}-\eqref{exp14} is bounded by
\begin{align*}
&C\, \norm{q}_{L^2}^2 \int_{\mathbb S^1} \int_0^{\infty}
h(s)^{2\b}\int_0^{\pi} \int_{-\pi}^{\pi}\int_0^{\pi}
\int_{Y^{\ast}(\Theta , \zeta' , \g' )} \sin\Theta \sin\zeta'\\
&\times |\q(B(j,s,\Theta , \t , \zeta' , \g'))|^2 d\t \, d\Theta
\, d\g' \, d\zeta' \, s\, ds \, d\s(z)\\
&= C\, \norm{q}_{L^2}^2 \int_0^{\infty} h(s)^{2\b}\int_0^{\pi}
\int_{-\pi}^{\pi}\int_0^{\pi}
\int_{Y^{\ast}(\Theta , \zeta' , \g' )} \sin\Theta \sin\zeta'\\
&\quad\times |\q(B(j,s,\Theta , \t , \zeta' , \g'))|^2 d\t \,
d\Theta \, d\g' \, d\zeta' \, s\, ds \, ,
\end{align*}
where
\begin{align*}
B(j,s,\Theta , \t , \zeta' , \g')&:={-h(s)\over 2}\, ((\sin\Theta
\,
\cos\t \, + \, \sin\zeta' \, \cos\g' )\, e_1\\
&\quad + \,(\sin\Theta \,\sin \t \, +\, \sin\zeta' \, \sin\g' )\,
e_2 \, + \, ( \cos\Theta \, +\, \cos\zeta' ) \, e_3).
\end{align*}
Next for each $j, \zeta' , \g'$ fixed we change variables $(s ,
\Theta , \t)\,\ra\, \mu = (\mu_1 , \mu_2 , \mu_3 )$ given by
\[
\mu = B(j,s,\Theta , \t , \zeta' , \g') .
\]
The  Jacobean of this transformation is given by
\begin{align*}
\left| { \partial (\mu_1 , \mu_2 , \mu_3 )\over \partial (s,
\Theta , \t ) } \right| &= {s \, h(s) \sin\Theta\over 8}\, ( 1 +
\sin\Theta\, \sin\zeta'\,\cos\t\, \cos\g'\\
&\quad + \, \sin\Theta\, \sin\zeta'\,\sin\t\, \sin\g'\, +\,
\cos\Theta\, \cos\zeta' )\, .
\end{align*}
Notice that the mentioned change involves an expression for $s$ in
terms of $\mu$ which depends on the parameters $j$, $\zeta'$,
$\g'$. Hence, the function $h$ has the same parametric dependence
\[
h(s)=h(\mu , j , \zeta' , \g' ) .
\]
Nevertheless, for $\t\in Y^{\ast}(\Theta , \zeta', \g')$ it
holds
\begin{equation}
\label{(12)}h(s)\sim |\mu|=|B(j,s,\Theta , \t , \zeta' , \g')|
\end{equation}
thanks to the condition stated in \eqref{(11a)}-\eqref{(11b)},
where
\begin{align*}
|B(j,s,\Theta , \t , \zeta' , \g')| &= {|h(s)|\over 2}\, ( 1 +
\sin\Theta\, \sin\zeta'\,\cos\t\, \cos\g'\\
&\quad + \, \sin\Theta\, \sin\zeta'\,\sin\t\, \sin\g'\, +\,
\cos\Theta\, \cos\zeta' )^{\um} .
\end{align*}
Indeed, property \eqref{(12)} is a reminiscence of the condition
${\e\over 100}\leq |\eta-\phi'-\xi'|\leq\e $ held for $\xi'\in
Y_{\eta}(\phi')$, see \eqref{3dexp48}.

We conclude that
\begin{align*}
G_j(\xi) &\leq C\,\norm{q}_{L^2}^2 \int_0^{\pi}\int_{-\pi}^{\pi
}\int_{\Rtres}|\mu|^{2\b-1}|\q(\mu)|^2 d\mu\, d\g'\, d\zeta'\\
&= C\,\norm{q}_{L^2}^2 \int_{\Rtres}|\mu|^{2\b-1}|\q(\mu)|^2
d\mu\, ,
\end{align*}
and then, we have proved the estimate
\begin{align*}
&\norm{R_{J_{j}}(q)}_{\homob}\leq C\norm{q}_{\homomenos
}\norm{q}_{L^2 }^2\norm{q}_{\dot{W}^{\b-{1\over 2}\, ,\, 2}}\,
,\\
\intertext{and since $\mathcal{N}$ is an universal constant we also have}\\
&\norm{Q_{I_2 }(q)}_{\homob}\leq C\norm{q}_{\homomenos
}\norm{q}_{L^2 }^2\norm{q}_{\dot{W}^{\b-{1\over 2}\, ,\, 2}}\, .
\end{align*}

This ends the proof of claim \eqref{3dexp45} and estimate
\eqref{afir2_3d}.

\qed

\proofof{estimate \eqref{afir3_3d} }This case is inspired on the
method used to control the piece $Q'_{II}(q)$ of the cubic term
from the Neumann-Born series in the three-dimensional case  in
\cite{RV}.

Let us start by decomposing the set $IV(\eta)$ as follows:
$IV(\eta)\subset IV_{<}(\eta)\cup IV_{>}(\eta)\,$, where
\begin{align*}
&IV_{<}(\eta):=\left\{ \tupla\in IV(\eta)\, :\, \xim , \taum ,
\phim\leq
\left( {1\over 50} + {1\over \sqrt{2}} \right)\e \right\}\, ,\\
&IV_{>}(\eta):=\left\{ \tupla\in IV(\eta)\, :\, |\eta-\xi|
,|\eta-\tau| , |\eta-\phi|\leq \left( {1\over 50} + {1\over
\sqrt{2}} \right)\e \right\}\, .
\end{align*}
In fact, if $\xim\leq {1\over\sqrt{2}}\,\e$ hence
\[
\phim\leq |\phi-\xi|+\xim\leq \left({1\over\sqrt{2}} + {1\over
100}\right)\e\, ,\quad \taum\leq |\tau-\phi|+\phim\leq \left(
{2\over 100} + {1\over\sqrt{2}}\right)\e\, ,
\]
and if $\xim\geq {1\over\sqrt{2}}\,\e$ then $|\eta-\xi|\leq
{1\over\sqrt{2}}\,\e\,$, and it holds
\[
|\eta-\phi|\leq |\eta-\xi|+|\xi-\phi |\leq \left({1\over\sqrt{2}}
+ {1\over 100}\right)\e\, ,\quad \etau\leq
|\eta-\phi|+|\phi-\tau|\leq \left( {2\over 100} +
{1\over\sqrt{2}}\right)\e\, .
\]

Taking the changes of variables $\xi=\eta-\xi'$,
$\tau=\eta-\tau'$, $\phi=\eta-\phi'$ in the integral
\[
\int\int\int_{IV_{>}(\eta)}
|\q(\xi')\q(\eta-\tau')\q(\tau'-\phi')\q(\phi'-\xi')|\,\measureee\,
,
\]
we notice that
$\wh{Q_{IV_{>}}(q)}(\eta)=\wh{Q_{IV_{<}}(q)}(\eta)\,$, and then $
\wh{Q_{IV}(q)}(\eta)\leq 2\, \wh{Q_{IV_{<}}(q)}(\eta)$. We take
another decomposition:
$IV_<(\eta)\subset\union_{k=1}^{+\infty}\left(IV_k(\eta)\cup
\wt{IV}_k(\eta) \right)\,$, being
\begin{align*}
&IV_k(\eta):=\left\{ \tupla\in IV_<(\eta)\, :\, \phim\leq \xim\sim 2^{-k}\e \right\}\, ,\\
&\wt{IV}_k(\eta):=\left\{ \tupla\in IV_<(\eta)\, :\, \xim\leq
\phim\sim 2^{-k}\e \right\}\, .
\end{align*}

It holds
\[
\norm{Q_{IV}(q)}_{\homob}\leq 2\,\suma_{k=1}^{+\infty}\left(
\norm{Q_{IV_k}(q)}_{\homob} + \norm{Q_{\wt{IV}_k}(q)}_{\homob}
\right)\, ,
\]
hence estimate \eqref{afir3_3d} follows from the following claims:
\addtocounter{equation}{-33}
\begin{subequations}
\begin{align} \label{3dexp54}&\norm{Q_{IV_k}(q)}_{\homob}\leq
C\, 2^{-\epsi k}\,
\norm{q}_{\homomenos}\norm{q}_{\dot{W}^{-\epsi\, ,\,
2}}\norm{q}_{L^2 }\norm{q}_{\homobmem}\, ,\\
\label{3dexp55a}&\norm{Q_{\wt{IV}_k}(q)}_{\homob}\leq C\,
2^{-\epsi k}\, \norm{q}_{\homomenos}\norm{q}_{\dot{W}^{-\epsi\,
,\,
2}}\norm{q}_{L^2 }\norm{q}_{\homobmem}\\
\label{3dexp55b}&\quad\quad\quad\quad\quad\quad\,\, + C\,
2^{-\epsi k}\, \norm{q}_{\dot{W}^{-\um-\epsi\, ,\,
2}}\norm{q}_{L^2}^2 \norm{q}_{\homobmem}\, ,
\end{align}\end{subequations}\addtocounter{equation}{32}
provided that $\epsi>0\,$.\\

\proofof{claim \eqref{3dexp54} }By Cauchy-Schwartz inequality,
\begin{align*}
\wh{Q_{IV_k}(q)}(\eta)&\leq {1\over\e^3}\,
\left(\int\int\int_{IV_k(\eta)}|\q(\xi)\q(\eta-\tau)|^2
\measure\right)^{\um}\\
&\quad\times
\left(\int\int\int_{IV_k(\eta)}|\q(\tau'-\phi')\q(\phi'-\xi')|^2
\measureee\right)^{\um} ,
\end{align*}
where
\begin{align*}
IV_k(\eta) &=\{\tupla\in\dominio : |\phi-\xi|\leq \e / 100 , \,
|\phi-\tau|\leq \e / 100 ,\\
&\quad\quad\xim , \taum , \phim\leq (1 / \sqrt{2} \, + \, 1/ 50
)\e , \, \phim\leq \xim\sim 2^{-k}\e\}.
\end{align*}
For each $\eta\in\Rtres\, ,\,\phi'\in\Gae$ fixed, using as above
the maximal operator we have
\[
\int_{\Gae}|\q(\tau'-\phi')|^2 d\s(\tau')\leq C\norm{q}_{L^2}^2\,
.
\]
Since $\s(\Gae)=\pi\,\e^2$ we also get
\begin{align}
\label{3dexp56}\norm{Q_{IV_k}(q)}_{\homob}^2&\leq C
\norm{q}_{L^2}^2\int_{\Rtres}\e^{2\b-4}\int_{B_k(\eta)}|\q(\xi)|^2
\int_{\mathcal{B}_{\xi}(\eta)}|\q(\eta-\tau)|^2
d\s(\tau)\, d\s(\xi)\\
\label{3dexp57}&\quad\times \int\int_{\mathcal{C}_k(\eta)}
|\q(\phi'-\xi')|^2
d\s(\xi') d\s(\phi') d\eta\\
\label{3dexp58}&\leq C \norm{q}_{\homomenos}^2
\norm{q}_{L^2}^2\int_{\Rtres}\e^{2\b-4}\int_{B_k(\eta)}|\q(\xi)|^2
\\
\label{3dexp59}&\quad\times 2^{-k}\e^2
\int_{\mathcal{B}_{\xi}(\eta)}|\q(\eta-\tau)|^2
d\s(\tau) d\s(\xi) d\eta\\
\label{3dexp60}&\leq C 2^{-2k\epsi}\norm{q}_{\homomenos}^2
\norm{q}_{L^2}^2 \int_{\Rtres}|\q(\xi)|^2
\int_{\Lambda^{\ast}(\xi)}|\eta|^{2\b-2}\\
\label{3dexp61}&\quad\times {\e^{2\epsi}\over\xim^{2\epsi}}
\int_{\mathcal{B}_{\xi}(\eta)}|\q(\eta-\tau)|^2 d\s(\tau)d\s(\eta) d\xi\\
\label{3dexp62}&\leq C 2^{-2k\epsi}\norm{q}_{\homomenos}^2
\norm{q}_{\dot{W}^{-\epsi\, ,\, 2}}^2
\norm{q}_{L^2}^2\norm{q}_{\homobmem}^2
\end{align}
where
\begin{align}
\notag&B_k(\eta):=\left\{ \xi\in\Gae\, :\, \xim\sim 2^{-k}\e\, ,\,
\xim\leq
\left( {1\over 50} + {1\over\sqrt{2}} \right)\e \right\} ,\\
\label{3dexp64}&\mathcal{B}_{\xi}(\eta):=\left\{ \tau\in\Gae\, :\,
\xitau\leq {\e\over 50}\, ,\, \taum\leq
\left({1\over 50} + {1\over\sqrt{2}}\right)\e \right\}\, ,\\
\notag&\mathcal{C}_k(\eta):=\left\{ (\xi',\phi')\in\Gae^2\, :\,
|\phi'-\xi'|\leq \metac\, ,\, |\phi'|\leq |\xi'|\sim 2^{-k}\e
\right\} ,\\
\label{exp1}&\Lambda^{\ast}(\xi):=\left\{\eta\in\Lax\, :\,
\xim\leq \left( {1\over 50} + {1\over\sqrt{2}} \right)\e \right\}
.
\end{align}
The estimate \eqref{f5} allows us to estimate
\eqref{3dexp56}-\eqref{3dexp57} by
\eqref{3dexp58}-\eqref{3dexp59}. The step from
\eqref{3dexp58}-\eqref{3dexp59} to \eqref{3dexp60}-\eqref{3dexp61}
follows from the property $\xim\sim 2^{-k}\e\,$ and the change of
the order of integration in $\xi$ and $\eta$ through Lemma
\ref{lema2.2}. Finally, we bound \eqref{3dexp60}-\eqref{3dexp61}
by \eqref{3dexp62} applying part $(ii)$ of Lemma \ref{lemmasum}.

\qed

\textbf{Proof of claim \eqref{3dexp55a}-\eqref{3dexp55b}.}

We split
$\wt{IV}_k(\eta)=\wt{IV}_{k,a}(\eta)\cup\wt{IV}_{k,b}(\eta)\, $,
where
\begin{align*}
&\wt{IV}_{k,a}(\eta):=\left\{ \tupla\in\wt{IV}_k(\eta)\, :\, \xim\geq 2^{-k-2}\e \right\}\\
&\quad\quad \quad\quad =\{\tupla\in\dominio : |\phi-\xi|\leq\e / 100 , \, |\phi-\tau|\leq \e / 100 ,\\
&\quad\quad\quad\quad\quad \quad \xim , \taum , \phim\leq
(1/\sqrt{2} \, +\, 1/ 50 )\e ,\, 2^{-k-2}\e\leq
\xim\leq \phim\sim 2^{-k}\e   \}\, ,\\
&\wt{IV}_{k,b}(\eta):=\left\{ \tupla\in\wt{IV}_k(\eta)\, :\,
\xim\leq 2^{-k-2}\e \right\}\\
&\quad\quad \quad\quad =\{\tupla\in\dominio : |\phi-\xi|\leq\e / 100 , \, |\phi-\tau|\leq \e / 100 ,\\
&\quad\quad\quad\quad\quad \quad \xim , \taum , \phim\leq
(1/\sqrt{2} \, +\, 1/ 50 )\e ,\, \xim\leq \phim\sim 2^{-k}\e , \,
\xim\leq 2^{-k-2}\e \}.
\end{align*}

On the domain $\wt{IV}_{k,a}(\eta)$, $2^{-k-2}\e\leq
\xim\leq\phim\leq 2^{-k+1}\e\,$ holds; hence following the steps
of the proof for the domain $IV_k(\eta)$ we arrive at
\begin{equation}
\label{3dexp63}\norm{Q_{\wt{IV}_{k,a}}(q)}_{\homob}\leq C\,
2^{-\epsi k}\, \norm{q}_{\homomenos}\norm{q}_{\dot{W}^{-\epsi\,
,\, 2}}\norm{q}_{L^2 }\norm{q}_{\homobmem}\, .
\end{equation}

On the domain $\wt{IV}_{k,b}(\eta)$ we have $|\xi-\phi|\geq
2^{-k-2}\e\,$, in fact
\[
|\xi-\phi|\geq \phim-\xim\geq 2^{-k-1}\e-2^{-k-2}\e = 2^{-k-2}\e\,
.
\]

In this case, we can bound
$\norm{Q_{\wt{IV}_{k,b}}(q)}_{\homob}^2$ by a similar expression
to \eqref{3dexp56}-\eqref{3dexp57} replacing $B_k(\eta)$ by the
set $\{\xi\in\Gae : \xim\leq 2^{-k-2}\e ,\, \xim \leq ({1\over 50}
+ {1\over \sqrt{2}})\e\}$ and the domain $\mathcal{C}_k(\eta)$ by
the set
\[
\left\{ (\xi',\phi')\in\Gae^2\, :\, 2^{-k-2}\e\leq
|\phi'-\xi'|\leq \metac\, ,\, |\xi'|\leq |\phi'|\sim 2^{-k}\e
\right\} .
\]
By the estimate \eqref{f6}, changing the order of integration in
$\xi$ and $\eta$ by Lemma \ref{lema2.2}, multiplying and dividing
by ${\e^{2\epsi}\over\xim^{2\epsi}}$ and applying that
${\xim^{2\epsi}\over\e^{2\epsi}}\leq 2^{2\epsi (-k-2)}\,$, and
finally by part $(ii)$ of Lemma \ref{lemmasum} we obtain
\begin{equation}
\label{3dexp65}\norm{Q_{\wt{IV}_{k,b}}(q)}_{\homob}\leq C \,
2^{-\epsi k}\, \norm{q}_{\dot{W}^{-\um-\epsi\, ,\,
2}}\norm{q}_{L^2}^2 \norm{q}_{\homobmem}\, .
\end{equation}

The expressions \eqref{3dexp63} and \eqref{3dexp65} lead up to
claim \eqref{3dexp55a}-\eqref{3dexp55b}.

This ends the proof of estimate \eqref{afir3_3d} and Proposition
\ref{propprincipal}.

\qed

\section{Proof of Theorem \ref{lestimate} (remainder term $ \a<n/2$).}

The control of the remainder term ${\bf R}_l$, where $l$ is as in
the statement of Theorem \ref{lestimate}, follows from the  next
proposition  by choosing  $C_0$ large enough in (\ref{qtilde}).
For $C_0 > (2\norm{q}_{\nohomoa})^4 $ we obtain the convergence of
the series
\[
\suma_{j=l}^{+\infty}\left( 2 C_0^{-{1\over
4}}\,\norm{q}_{\nohomoa}\right)^j\, .
\]
\begin{propo}
\label{jestimate}{\it Let $n\in\{2,3\}$, $q\in W^{\a , 2}(\Rn)$
compactly supported and $0\leq \a < n/2$. Assume that $C_0
>1\,$, $j\geq 4$ if $n=2$ and $j\geq 5$ if $n=3$. Then, for any $\b\in\R$ such that $\b<\a+1$, it holds:
\begin{equation}
\label{triang}\norm{\widetilde{Q}_j (q)}_{\homob}\leq C(\a ,\b )\,
C_0^{17 /4} (2C_0^{-{1\over 4}}\norm{q}_{\nohomoa})^j
\,\norm{q}_{L^2 }\norm{q}_{\nohomoa}^{-1} .
\end{equation}
 }
\end{propo}

\proofof{Proposition \ref{jestimate}} We write
$R_{\t,k}f(x)=e^{-ik\t\cdot x}R_k(e^{ik\t\cdot
(\cdot)}f(\cdot))(x)$. It holds
\[
\widehat{Q_{j}(q)} (\xi)=\int_{\Rn}e^{ik\t \cdot
y}(qR_{k})^{j-1}(q(\cdot)e^{ik\t \cdot (\cdot)})(y)dy =
\int_{\Rn}e^{2ik\t\cdot y}(q R_{\t , k})^{j-1}(q)(y)\, dy \, ,
\]
where $k=\xim /2$ and $\t=-\xi / \xim$. Let $\g\in\R$ be such that
$\g< \b_j$, being
\[
\b_{j}:=\left\{
\begin{array}{llll}
{3\over 4}(j-2)+{\a\over 4}(j-1)\, , & \text{if }\a\leq {1\over
2}\,\text{ and }\,\, n=2\, ,\\
(j-3)({3\over 4}+{\a\over 4})+1\, , & \text{if }{1\over 2}\leq
\a\leq 1\,\text{ and }\,\, n=2\, ,\\
{j-2\over 2}+(j-1){\a\over 3}-\um\, , &\text{if }0\leq \a\leq
{3\over 4}\,\text{ and }\,\, n=3\, ,\\
(j-3)\left(\um+{\a\over 3}\right)+\um\, , &\text{if }{3\over
4}\leq \a\leq {3\over 2} \,\text{ and }\,\, n=3\, .
\end{array}
\right.
\]
Taking the change of variables $\xi=-2k\t\,$ with $k\geq 0 \,$,
$\t\in S^{n-1}$, we have
\begin{align}
\notag\|\widetilde{Q}_{j}(q)\|_{\homog}^2
&=\int_{\Rn}\xim^{2\g}|\widehat{\widetilde{Q}_{j}(q)}(\xi)|^2 d\xi
= \int_{\Rn}\xim^{2\g}\chi^{\ast}(\xi)
|\widehat{Q_{j}(q)}(\xi)|^2 d\xi\\
\label{2dexp14}&\leq C_n\, 2^{2\g}\int_{k={C_{0}\over
2}}^{+\infty}k^{2\g+n-1}\int_{S^{n-1} }\|(q R_{\t ,
k})^{j-1}(q)\|_{L^{1}}^2 d\s(\t)dk\, .
\end{align}
By Cauchy-Schwartz inequality, $ \|(q R_{\t ,
k})^{j-1}(q)\|_{L^{1}} \leq C \|q\|_{L^{2}}\|R_{\t , k}(q R_{\t ,
k})^{j-2}(q)\|_{L^{2}}$. Using the estimate given by Lemma 3.4 in
\cite{R1} for the operator $\rteta$ and the next inequality for
products of Sobolev spaces due to Zolesio
\[
\norm{fg}_{W^{\a_3 , p}}\leq C\, \norm{f}_{W^{\a_1 ,
p_1}}\norm{g}_{W^{\a_2 , p_2}} ,
\]
where $\a_1 ,\a_2 , \a_3\geq 0$, $\a_3\leq \a_j$, $p_j > p,$
$j=1,2$, $\a_1 + \a_2 - \a_3 \geq n (1/p_1 + 1/ p_2 - 1/p )\geq
0$, we   arrive at
\[
\|\rteta(q\rteta)^{j-2}(q)\|_{L^{2}}\leq C
k^{\g_j}\|q\|_{\nohomoa}^{j-1}\, ,
\]
where $\g_j:=-(j-1)+{n-1\over 2}(j-3)(1/2-\a/n)+{n-1\over
2}\max\{0,\, {1\over 2}-{2\a\over n} \}\, $. All this leads us up
to
\begin{equation*}
\|\widetilde{Q}_{j}(q)\|_{\homog}^2 \leq C \,
2^{2\g}\int_{k={C_{0}\over 2}}^{+\infty}k^{2\g +n-1+2\g_j}dk\,
\|q\|_{L^{2}}^2 \|q\|_{\nohomoa}^{2j-2}\, ,
\end{equation*}
where the integral converges if $2\g+2\g_j+n<0$, that is to say,
if $\g < \b_j$. Notice that $\b_j=-{n\over 2}-\g_j$. In this way,
we have proved that, for $\g<\b_j$, it holds
\begin{equation}
\label{2dexp15}\|\widetilde{Q}_{j}(q)\|_{\homog} \leq C
 {2^{\g}\over
\sqrt{\b_j-\g}} \left( \f{C_{0}}{2} \right)^{\g-\b_j} \,
\|q\|_{L^{2}} \|q\|_{\nohomoa}^{j-1}\, .
\end{equation}

Let $\epsi=\epsi (\a , \b ):= (\a +1)-\b
>0\,$. Keeping in mind $2\b=2(\b_j-\epsi)+2(\a+1-\b_j )\,$ and  $\a+1\leq \b_j\,$ for our $j$, we write
\begin{align*}
\norm{\widetilde{Q}_j (q )}_{\homob} &\leq C_0 ^{\a +1-\b_j }
\norm{\widetilde{Q}_j
(q)}_{\dot{W}^{\b_j -\epsi , 2}} \\
&\leq C\, C_0 ^{\a +1-\b_j } {2^{\b_j-\epsi}\over
\sqrt{\epsi}}\,\left({C_0\over 2}\right)^{-\epsi} \norm{q}_{L^2 }
\norm{q}_{\nohomoa}^{j-1}\\
&= C(\a , \b ) \, 2^{\b_j } \, C_0^{\b -\b_j } \norm{q}_{L^2 }
\norm{q}_{\nohomoa}^{j-1}\, ,
\end{align*}
where the last inequality follows from formula \eqref{2dexp15} in
the case $\g = \b_j -\epsi\,$. In our setting, $2^{\b_j } \leq
2^{j }\,$. Moreover, $\b -\b_j < \a+1 - \b_j \leq -{1\over
4}j+{17\over 4}\,$, for our $j$, and $C_0>1$. We have proved
\eqref{triang}.

\qed

\section{The case  $\a\geq n/2$.}

In this section we are going to extend Theorem \ref{theoremone},
Theorem \ref{lestimate} and estimates \eqref{A}, \eqref{B},
\eqref{C} and \eqref{D1}-\eqref{D2} for any $\a\geq 0$. Then
Theorem \ref{theorem3} will follow from these estimates for any
$\a\geq 0$. We start with a Leibniz' type formula for derivatives
of $Q_j(q)$ which we state as follows
\begin{teorema}
\label{teor_partes}{\it Assume that $\a\in\mathbb{N}^n \, ,\,
j\in\Z\, ,\, j\geq 2\,$ and let $q\in W^{|\a| , 2}(\Rn)$ be a
compactly supported function. Then
\[
D^{\a}Q_j(q)= \suma_{\substack{\b_1 + \,\dots\, + \b_{j}=\a\\
\b_1 ,\, ... \, , \b_j\geq 0}} { \a! \over
\b_1!\cdot\,\dots\,\cdot \b_{j}! }\, Q_j(D^{\b_1}q \, ,\, \dots \,
,\, D^{\b_{j}}q)\, .
\]

}
\end{teorema}

\begin{nota}
From the proof of Theorem \ref{teor_partes} one also deduces the
formula
\begin{equation}
\label{apform9}D^{\a}\wt{Q}_j(q)= \suma_{\substack{\b_1 +
\,\dots\, + \b_{j}=\a\\ \b_1 ,\, ... \, \b_j\geq 0}} { \a! \over
\b_1!\cdot\,\dots\,\cdot \b_{j}! }\, \wt{Q}_j(D^{\b_1}q \, ,\,
\dots \, ,\, D^{\b_{j}}q)\, ,
\end{equation}
for the same hypotheses on $q$.
\end{nota}

\proofof{Theorem \ref{teor_partes}} Writing the resolvent $R_k$ as
the convolution operator with the outgoing fundamental solution to
the Helmholtz equation (see \cite{CK}, \cite{R2})
\begin{equation}
\label{solfund}\phi_k(x)=C_n
k^{(n-2)/2}{H^{(1)}_{(n-2)/2}(k|x|)\over |x|^{(n-2)/2} }\, ,
\end{equation}
where $C_n = {1\over 2i(2\pi)^{(n-2)/2}}\,$ and
$H^{(1)}_{(n-2)/2}$ denotes the Hankel function of the first kind
and order $(n-2)/2$ (see \cite{W}), we have
\begin{align*}
\wh{Q_j(q)}(-2k\t) &=\int_{\Rn} e^{ik\t\cdot y}q(y)( R_k q
)^{j-1} ( e^{ik \t\cdot (\cdot ) })(y) \, dy\\
&=\int_{\left( \Rn\right)^{j}} e^{ik\t\cdot \, x_1 }q(x_1 )\,
\displaystyle\prod_{l=1}^{j-1} \left(
\phi_k(x_l-x_{l+1})q(x_{l+1})\,\right)\, e^{ik\t\cdot\, x_{j}}\,
dx\, ,
\end{align*}
and
\begin{align}
\label{h1}&\mathcal{F}\left( Q_j (f_1 , \dots , f_{j}) \right) (-2k\t)\\
\label{h2}&=  \int_{\left(\Rn\right)^{j}} e^{ik\t\cdot \, x_1
}f_1(x_1 )\, \displaystyle\prod_{l=1}^{j-1} \left(
\phi_k(x_l-x_{l+1})f_{l+1}(x_{l+1})\,\right)\, e^{ik\t\cdot\,
x_{j}}\, dx\, ,
\end{align}
where $dx=dx_1 \cdot\, \dots \, \cdot  dx_{j}\,$ and
$x_l\in\Rn\,$, for any $l=1,...,j\,$. We know that
$\mathcal{F}\left( D^{\a}Q_j(q)\right)(-2k\t)=(-i 2 k \t
)^{\a}\wh{Q_j(q)} (-2k\t )\,$. Taking the change $x_l=x_1 + y_l$,
$2\leq l\leq j$, it holds
\begin{align*}
&\wh{Q_j(q)} (-2k\t )\\
&= \int_{\left(\Rn\right)^{j}} e^{i2k\t\cdot \, x_1 }q(x_1 )\,
\displaystyle\prod_{l=1}^{j-1} \left( \phi_k(x_l
-x_{l+1})q(x_{l+1}) e^{-ik\t\cdot\, (x_l
-x_{l+1})}\,\right)\, dx\\
&= \int_{\Rn}\int_{\left(\Rn\right)^{j-1}} e^{i2k\t \cdot \, x_1
}q(x_1 )\,\displaystyle\prod_{l=1}^{j-1} \left( \phi_k(y_l
-y_{l+1})q(x_1+y_{l+1}) e^{-ik\t\cdot\, (y_l -y_{l+1})}\,\right)\,
dy dx_1 \, ,
\end{align*}
where $y_1 =0$ and $dy = dy_j\cdot\,\dots\,\cdot dy_2$.
Integrating by parts, we have
\begin{align*}
&\left(-i2k\t \right)^{\a} \wh{Q_j(q)} (-2k\t)\\
&= (-1)^{|\a|}\int_{\Rn}\int_{\left(\Rn\right)^{j-1} }
(D^{\a}_{x_1 } e^{i2k\t\cdot \, x_1 })\, q(x_1
)\,\displaystyle\prod_{l=1}^{j-1} \left( \phi_k(y_l
-y_{l+1})q(x_1+y_{l+1}) e^{-ik\t\cdot\, (y_l -y_{l+1})}\,\right)\,
dy dx_1\\
&= \int_{\Rn}\int_{\left(\Rn\right)^{j-1} } e^{i2k\t\cdot \, x_1
}\, D^{\a}_{x_1} \left[ q(x_1 )\,\displaystyle\prod_{l=1}^{j-1}
\left( \phi_k(y_l -y_{l+1})q(x_1+y_{l+1}) e^{-ik\t\cdot\, (y_l
-y_{l+1})}\,\right)\right]\, dy dx_1\\
&= \suma_{\substack{\b_1 + \,\dots\, + \b_{j}=\a\\
\b_1 ,\, ... \, , \b_j\geq 0}} { \a! \over
\b_1!\cdot\,\dots\,\cdot \b_{j}! }
\int_{\Rn}\int_{\left(\Rn\right)^{j-1} }e^{i2k\t\cdot \, x_1 } \,
D^{\b_1 } q(x_1 )\\
&\qquad\times \displaystyle\prod_{l=1}^{j-1} \left( \phi_k(y_l
-y_{l+1}) e^{-ik\t\cdot\, (y_l -y_{l+1})} D^{\b_{l+1}}
q(x_1+y_{l+1})\right)\, dy dx_1\, ,
\end{align*}
where we have applied Leibniz' formula:
\[
D^{\a}(f_1\cdot\, \dots\,\cdot f_k )=\suma_{\substack{\b_1 + \dots + \b_k =\a\\
\b_1 ,\, ... \, , \b_k\geq 0}}{\a!\over\b_1!\cdot\, \dots
\,\cdot\b_k! }\, D^{\b_1}f_1\cdot\, \dots \, \cdot D^{\b_k}f_k\, .
\]
Finally,
\begin{align*}
\left(-i2k\t \right)^{\a} \wh{Q_j(q)} (-2k\t)&=\suma_{\substack{\b_1 + \,\dots\, + \b_{j}=\a\\
\b_1 ,\, ... \, , \b_j\geq 0}} { \a! \over
\b_1!\cdot\,\dots\,\cdot \b_{j}! } \int_{\left(\Rn\right)^{j}}
e^{ik\t\cdot \, x_1 }D^{\b_1}q(x_1
)\\
&\qquad\times\displaystyle\prod_{l=1}^{j-1} \left( \phi_k(x_l
-x_{l+1})D^{\b_{l+1}}q(x_{l+1}) \,\right) e^{ik\t\cdot\, x_{j}}\,
dx\, ,
\end{align*}
and remembering the expression \eqref{h1}-\eqref{h2}, we end the
proof of Theorem \ref{teor_partes}.

\qed

To prove Theorem \ref{theorem3} in case $\a\geq n/2$  we use
induction on   $[\a]$. We need to use  the boundedness  of  the
$j$-multiple scattering operators, see the notation, when acting
on $(q_1,...q_j)$ where $j-1$ of them are equal to $q$. Namely
\begin{propo}
\label{teo5}{\it

Let $n\in\{2,3\}$, $\a\in\R\,$ with $0\leq \a <n / 2\,$ and let us
suppose that $q_1 , q_2 \in W^{\a\, , \, 2}(\Rn)$ are compactly
supported functions. Then $Q_2(q_1 , q_2 )\in W^{\b\, ,\,
2}(\Rn)+\reg\,$, for any $\b\in\R$ such that $0\leq \b<\a +
{1\over 2}\,$. Moreover, there exists a constant $C(\a , \b , q_1
, q_2)>0$ which just depends of $\a,\b $ and the supports of $q_1
, q_2$ such that
\[
\norm{\wt{Q}_2(q_1 , q_2 )}_{\homob}\leq C(\a , \b , q_1 , q_2 )
\, \maxim\{\norm{q_1}_{W^{\a , 2}}^2 , \norm{q_2}_{W^{\a , 2}}^2
\}\, .
\]

}
\end{propo}

Proposition \ref{teo5} follows by polarization of estimate
\eqref{A}.

Next propositions \ref{teo6_3d}, \ref{teo7} and
\ref{jestimatemulti} follow from the proofs of the analogous
estimates \eqref{B}, \eqref{C}, \eqref{D1}-\eqref{D2}  and
Proposition \ref{jestimate}.

\begin{propo}
\label{teo6_3d}{\it

Let $n\in\{2,3\}$, $\a\in\R\,$ with $0\leq \a <n/2\,$ and $q_1 ,
q_2 , q_3 $ as $q_1$ from Proposition \ref{teo5}. Then $Q_3(q_1 ,
q_2 , q_3 )\in W^{\b\, ,\, 2}(\Rn) + C^{\infty }(\Rn)\,$, for any
$\b\in\R$ holding $0\leq\b <\a +1$ if $n=2$ and $0\leq \b<\a +
1/2\,$ if $n=3$. Moreover, there exists a constant $C(\a , \b ,
q_1 , q_2 , q_3 )$ that just depends of $\a,\b $ and $\supp q_1 ,
\supp q_2 , \supp q_3$ such that
\begin{align}
&\norm{\wt{Q}_3(q_1 , q_2 , q_3 )}_{\homob}\\
\label{exp5}&\leq C(\a , \b , q_1 , q_2 , q_3 )\, (
 \suma_{\s\in S_3} \|q_{\s (1)}\|_{L^2 }\|q_{\s (2)}\|_{L^2}\|q_{\s (3)}\|_{\homoa
 }\\
\label{exp6}&\quad\quad\,\,\,\quad\quad\quad\quad\quad +
\suma_{\tau\in S_3} \|q_{\tau (1)}\|_{\homomenos }\|q_{\tau
(2)}\|_{\dot{W}^{-\epsi ,
2}}\|q_{\tau (3)}\|_{\homoa }\\
\label{exp7}&\quad\quad\,\,\,\quad\quad\quad\quad\quad +
 \suma_{\o\in S_3} \|q_{\o (1)}\|_{\dot{W}^{-\um-\epsi , 2} }\|q_{\o (2)}\|_{L^2}\|q_{\o (3)}\|_{\homoa
 }\, )\, ,
\end{align}
where $\epsi = \a +1-\b>0$ if $n=2$ and $\epsi=\a+\um-\b>0$ if
$n=3$.

}
\end{propo}

\begin{propo}
\label{teo7}{\it

Let $\a\in\R\,$ with $0\leq \a <3/2\,$ and $q_1 , q_2 , q_3 , q_4$
as $q_1$ from Proposition \ref{teo5} for $n=3$. Then $Q_4(q_1 ,
q_2 , q_3 ,q_4)\in W^{\b\, ,\, 2}(\Rtres) + C^{\infty
}(\Rtres)\,$, for any $\b\in\R$ with $0\leq \b<\a +1/2\,$.
Moreover, there exists a constant $C(\a , \b , q_1 , q_2 , q_3 ,
q_4 )>0$ just depending of $\a ,\b $ and the supports of $q_1 ,
q_2 , q_3 , q_4$ such that
\begin{align}
\notag&\norm{\wt{Q}_4(q_1 , q_2 , q_3 , q_4 )}_{\homob}\\
\label{exp8}&\leq C(\a , \b , q_1 , q_2 , q_3 , q_4 )\, (
 \suma_{\s\in S_4} \|q_{\s (1)}\|_{L^2 }\|q_{\s (2)}\|_{L^2}\|q_{\s
 (3)}\|_{L^2
 }\|q_{\s (4)}\|_{\homoa
 }\\
\label{exp9}&\quad\,\,\quad\quad\,\,\,\quad\quad\quad\quad\quad +
\suma_{\tau\in S_4} \|q_{\tau (1)}\|_{\homomenos }\|q_{\tau
(2)}\|_{L^2}\|q_{\tau (3)}\|_{L^2 }\|q_{\tau (4)}\|_{\homoa
 }\\
\label{exp10}&\quad\,\,\quad\quad\,\,\,\quad\quad\quad\quad\quad +
 \suma_{\o\in S_4} \|q_{\o (1)}\|_{\homomenos }\|q_{\o (2)}\|_{\dot{W}^{-\epsi , 2}}\|q_{\o
 (3)}\|_{L^2 }\|q_{\o (4)}\|_{\homoa
 }\\
 \label{exp11}&\quad\,\,\quad\quad\,\,\,\quad\quad\quad\quad\quad +
 \suma_{\rho\in S_4} \|q_{\rho (1)}\|_{\dot{W}^{-\um-\epsi , 2} }\|q_{\rho (2)}\|_{L^2}\|q_{\rho
 (3)}\|_{L^2
 }\|q_{\rho (4)}\|_{\homoa
 }\, )\, ,
\end{align}
where $\epsi=\a+\um-\b>0$.

}
\end{propo}

\begin{propo}
\label{jestimatemulti} {\it Let us assume that $n\in\{2,3\}\,$,
$\a\in\R$, $0\leq \a <n/2\,$, $q_1 , ... , q_j \in W^{\a\, ,\,
2}(\Rn)$ are compactly supported functions and $C_0>1$. Hence, for
any $\b\in\R$ such that $\b<\a+1\,$:
\begin{equation}
\label{exp4}\norm{\widetilde{Q}_j (q_1 , ... , q_j)}_{\homob}\leq
C(\a , \b )\,C_0^{17/4}\left( 2 C_0^{-{1\over 4}}\,\maxim_{1\leq
l\leq j}\norm{q_l}_{\nohomoa}\right)^j\, ,
\end{equation}
where $j\geq 4\,$ if $n=2$ and $j\geq 5\,$ if $n=3$. }

\end{propo}

\section{Appendix.}

In this section we state two results, Lemma \ref{lema2.2} and
\ref{remark}, which are often used in this work and state and
prove an important result, Lemma \ref{lemmasum}, in order to
demonstrate Proposition \ref{propprincipal}.

Let $V$ be the  submanifold of $\mathbb{R}^{2n}\,$ $V:=\{
(\eta,\xi) \in \Rn \times \Rn : \xi\cdot (\xi-\eta)=0 \}.$  Then
$V$ can be
 viewed as a bundle of  spherical
sections $ V=\{ (\eta,\xi) \in \Rn \times \Rn : \eta \in \Rn ,
\xi\in\Gae \}$, or as  a bundle of  plane sections: $V=\{
(\eta,\xi) \in \Rn \times \Rn : \xi \in \Rn, \eta\in\Lax \}\,$,
where $\Gae$ and $\Lax$ are defined in \eqref{notation1} and
\eqref{notation2}. In this context, the following lemma from
\cite{RV} allows us to change the order of integration in $\xi$
and $\eta\,$.
\begin{lema}
\label{lema2.2}{\it Let $V= \{ (\eta,\xi) \in \Rn \times \Rn :
\xi\cdot (\xi-\eta)=0 \}$. Let $d\s_{\eta}(\xi)$ be the measure on
$\Gae$ induced by the $n$-dimensional Lebesgue measure $d\xi$ and
let $d\s_{\xi}(\eta)$ be the measure on $\Lambda(\xi)$ induced by
the $n$-dimensional Lebesgue measure $d\eta$. Then
\[
d\s_{\eta}(\xi)d\eta={\e\over \xim}\, d\s_{\xi}(\eta)d\xi.
\]}
\end{lema}

The following lemma in \cite{RV} is used several times in this
work.
\begin{lema}
\label{remark}{\it Assume that the support of $q$ is contained in
the unit ball. Then we have:

$(1)$ If $\,\xi, \ \xi' \in \Rn$ satisfy $|\xi-\xi'|\leq 3$, then
$ |\q(\xi)|\leq CM\q(\xi')\,$.

$(2)$ $\| \q \|_{L^{\infty}}\leq C\| \q\|_{L^{2}}$.

$(3)$ For $0<\g<{n\over 2}$, $ \| q \|_{\dot W^{-\g,\, 2}}\leq C
\| q \|_{L^{2}}\, $, where $C$ depends
 on the size of the support of $q$.}
\end{lema}

The following lemma is fundamental to control the spherical term
$Q (q)$ of the quartic term.

\begin{lema}
\label{lemmasum} {\it Let $\xi\in\Rtres\setminus\{0\},\, \b
\in\R\, , \, \epsi>0 ,\, k\in\N$. We denote
\begin{align}
\label{3dexp52}&F_k(\xi):=\int_{\Lax}
\e^{2\b-4}\int\int_{A_k(\eta)}|\q(\eta-\phi'-\xi')|^2 d\s(\xi')
d\s(\phi') d\s(\eta)\, ,\\
\label{3dexp20}&G(\xi):=\int_{\Lambda^{\ast}(\xi)}\e^{2\b-2+2\epsi}\int_
{\mathcal{B}_{\xi}(\eta)}|\q(\eta-\tau)|^2 d\s(\tau) d\s(\eta)\, ,
\end{align}
where $A_{k}(\eta)\, $, $\mathcal{B}_{\xi}(\eta)$,
$\Lambda^{\ast}(\xi)$ are defined in \eqref{3dexp9},
\eqref{3dexp64}, \eqref{exp1} respectively. Then

$(i)$ $ F_k(\xi)\leq C\, 2^{-2k}
\int_{\Rtres}|\l|^{2\b-1}|\q(\l)|^2 d\l \, $, for some constant
$C$ independent of $\xi$ , $k$ and  $q$.

$(ii)$ $G(\xi)\leq C\int_{\Rtres}|\l|^{2\b-1+2\epsi}|\q(\l)|^2
d\l\,$, for some constant $C$ independent of $\xi$ and $q$.

}
\end{lema}

\proofof{Lemma \ref{lemmasum}}

$\bullet$ Proof of $(i)$. For $\eta\in\Lax $ we write $\eta=\xi +s
z\,$ and $h(s):=\e = \left( \xim^2 +s^2 \right)^{{1\over 2}}$,
where $s\geq 0$ and $z\in\{\xi\}^{\bot}\,$, $|z|=1\,$. For
simplicity, we don't specify the dependence of variables with
respect to $\xi$ since it is fixed along the proof  . It holds
$d\s(\eta)=s\, ds d\s(z)\,$, where $d\s(z)$ denotes the measure on
the unitary circumference $\mathbb S^1$ in the plane
$\{\xi\}^{\bot}\,$. We have
\begin{equation}
\label{3dexp4}F_k(\xi)=\int_0^{\infty }\int_{\mathbb S^1
}h(s)^{2\b-4 }\int\int_{A_k(s , z)}|\q(\xi+sz-\phi'-\xi')|^2
d\s(\xi')d\s(\phi')d\s(z)s\, ds\, .
\end{equation}
Fixing $z  , s$ we parametrize $\xi', \phi'\in \Gamma(\xi + s z)$
by $v , u\in \mathbb S^2$, respectively:
\[
\xi'={\xi +  s z\over 2} + {h(s)\over 2}\, v\, , \,\, \phi'={\xi +
s z\over 2} + {h(s)\over 2 }\, u\, ,\quad u , v\in \mathbb S^2 \,
,
\]
where $d\s(\xi')=C\, h(s)^2 d\s(v)\, ,\, d\s(\phi')=C\, h(s)^2
d\s(u)\,$. The domain of integration for $v, u$ is given by
\[
\{ (v , u) \in \mathbb S^2\times \mathbb S^2 \, :\, |u-v|\leq
2^{-k+2}\, ,\, 1+u\cdot v\geq 1/ 5000  \}\, ,
\]
since $|\xi'-\phi'|\leq 2^{-k+1}h(s)$ implies $|u-v|\leq
2^{-k+2}\,$ and $|\xi + s z-\phi'-\xi'|\geq {h(s)\over 100}$
implies $1+u\cdot v\geq {1\over 5000}\,$. Let $\{D_{j,k} :
j\in\{1,\, ...\, , N_0 2^{2k}\}\,\}$ be a finite overlapping cover
of the sphere $\mathbb S^2 $ with overlapping constant independent
of $k$ such that $D_{j,k}$ is an spherical cup of diameter
${2^{-k}\over 50}$ and $N_0$ an appropriate constant. For each $j$ we
define
\[
\wt{D}_{j,k} :=\left\{ u\in \mathbb S^2 : |u-v|\leq 2^{-k+2}\, ,\,
1+u\cdot v\geq {1\over 5000}\, ,\, \text{for some }v\in D_{j,k}
\right\}\, .
\]
The expression \eqref{3dexp4} is bounded by
\begin{equation}
\label{3dexp5}C \suma_{j=1}^{N_02^{2k}}\int_{\mathbb S^1
}\int_{D_{j,k} }\int_0^{+\infty}\int_{\wt{D}_{j,k}
}h(s)^{2\b}\left|\q(-{h(s)\over 2}\, \left(u+v\right) )\right|^2
d\s(u)\, s\, ds\, d\s(v)d\s(z)\, .
\end{equation}
Notice that for every $j\in\{1,\,\dots\, , N_02^{2k}\}$,
$u\in\wt{D}_{j,k}$ and $v\in D_{j,k}$, we have   $|u+v|\geq
{1\over 100}$. In fact, since $u\in\wt{D}_{j,k}$ there
exists a $v'\in D_{j,k}$ such that $1+u\cdot v'\geq 1/5000\,$, and
hence $|u+v'|\geq 1/50$. We have
\begin{align*}
|u+v|&\geq |u+v'| - |v-v'|\geq {1\over 50}-\text{diam}\,D_{j,k}\\
&={1\over 50} - {2^{-k}\over 50}\geq {1\over 50}-{1\over
100}={1\over 100}.
\end{align*}
We take spherical coordinates for $u$ with respect to the
canonical reference of $\Rtres$,
\begin{equation}
\label{esfericas}u=(\cos\t \sin\psi\, ,\, \sin\t\sin\psi\, ,\,
\cos\psi)\, ,
\end{equation}
with $d\s(u)=\sin\psi \, d\psi d\t\,$. We bound \eqref{3dexp5} by
\[
C \suma_{j=1}^{N_02^{2k}}\int_{\mathbb S^1 }\int_{D_{j,k}
}\int_0^{+\infty}\int\int_{D^{\ast}_{j,k}
}h(s)^{2\b}\left|\q(-{h(s)\over 2}\, \left(u(\psi , \t )+v\right)
)\right|^2 \sin\psi\, d\psi d\t \, s\, ds\, d\s(v)d\s(z)\, ,
\]
where $u(\psi , \t)$ is given by \eqref{esfericas} and
$D^{\ast}_{j,k} :=\{(\psi , \t )\in [0,\pi ] \times [0,2\pi )\, :
\, u(\psi , \t )\in \wt{D}_{j,k} \}\,$. We remark that for $(\psi
, \t )\in D^{\ast}_{j,k}$ it holds
\begin{equation}
\label{f10}|u(\psi , \t) + v|\geq {1\over 100}.
\end{equation}
For $z\in \mathbb S^1\, ,\, j\in\{1, ... , N_02^{2k}\}$ and $v\in
D_{j,k} $ fixed we consider the change $(s , \t , \psi
)\rightarrow (\l_1 , \l_2 , \l_3 )=\l$ given by
\begin{align*}
\l&=\xi + s z-\phi'-\xi'=-{h(s)\over 2}\, (u(\psi , \t )+v
)\\
&=-{h(s)\over 2}\, \left( \,\cos\t\sin\psi +v_1\, ,\,
\sin\t\sin\psi +v_2 \, ,\, \cos\psi +v_3  \, \right)\, ,
\end{align*}
where $v=(v_1 , v_2 , v_3)$. We have   $d\l = {sh(s)\sin\psi
|1+u(\psi , \t ) \cdot v |\over 8}\, dsd\psi d\t$. From
condition \eqref{f10} we deduce that $|\l |\sim h(s)\,$. We define
this family of sets with overlapping constant independent of $k$
contained in the interior of convex cones in $\Rtres$:
\[
H_{j,k}:=\{r(u+v )\, :\, u\in\wt{D}_{j,k} \, ,\, v \in D_{j,k} \,
, \, r<0 \}\, .
\]
Since $\s(D_{j,k} )\sim 2^{-2k}\,$, we have
\begin{align*}
F_k(\xi)&\leq C \suma_{j=1}^{N_02^{2k}}\int_{\mathbb S^1
}\int_{D_{j,k}
}\int_{H_{j,k} } |\l|^{2\b-1 } |\q(\l )|^2 d\l \, d\s(v)d\s(z)\\
&\leq C \, 2^{-2k } \int_{\Rtres } \left( \suma_{j=1}^{N_02^{2k}}
\chi_{H_{j,k} } (\l) \right) |\l|^{2\b-1 } |\q(\l )|^2 d\l\\
&\leq C \, 2^{-2k }\int_{\Rtres } |\l|^{2\b-1 }|\q(\l )|^2 d\l\, .
\end{align*}

\qed

$\bullet$ Proof of $(ii)$. We also express
$\eta\in\Lambda^{\ast}(\xi)$ as $\eta = \xi + s z\,$, with $s>0$,
$z\in\{\xi\}^{\bot}\,$, $|z|=1\,$. We use the same notation
$h(s)$. Since $\xim\leq\left(\cte\right) \e\,$, it is true that
$s\geq 0.9\xim\,$. It holds $d\s(\eta)=s ds d\s(z)\,$. For $z\in
\mathbb S^1$ fixed, we take a reference $\{e_1 , e_2 , e_3 \}$ in
$\Rtres$ such that $z=e_3\,$. Then $\xi = (\xi_1 , \xi_2 , 0)$ and
$\eta=(\xi_1 , \xi_2 , s)\,$. We write for $\tau\in
\mathcal{B}_{\xi}( \xi + s z)$,
\[
\tau={\xi + s z\over 2} + {h(s)\over 2} \, v\, , \quad v\in
\mathbb S^2\, ,
\]
with $d\s(\tau)=C\, h(s)^2 d\s(v)\,$. We express $v$ in spherical
coordinates with respect of our reference:
\[
v=(\,\sin\Theta\,\cos\,\t , \,\sin\Theta\,\sin\t , \,\cos\,\Theta
)\, , \quad 0\leq \Theta\leq \pi\, ,\,\, -\pi\leq \t < \pi\, ,
\]
where $d\s(v) = \sin \Theta\, d\Theta \, d\t$. For $s,z$ fixed,
let
\[
\mathcal{M}(s,z):=\{(\Theta , \t)\in [0, \pi]\times [-\pi , \pi)\,
:\, \tau(s,z,\Theta ,\t)\in \mathcal{B}_{\xi}(\xi + sz)\}\, .
\]
We obtain
\begin{align*}
G(\xi) &\leq \int_{\mathbb S^1 }\int_{0.9\xim
}^{\infty}\int\int_{\mathcal{M}(s,z)} s\, h(s)^{2\b
+2\epsi}\,\sin\Theta
\\
&\quad\times \left|\q\left( {1\over 2}\, \left(
\xi_1-h(s)\,\cos\,\t \,\sin\Theta , \xi_2 -h(s)\,\sin\t
\,\sin\Theta , s-h(s)\,\cos\,\Theta \right)\, \right)\right|^2
d\Theta d\t ds d\s(z)\, .
\end{align*}
For each $u$, we take the change of variables $(s,\Theta ,\t )\,
\rightarrow\, \l=(\l_1 , \l_2 , \l_3)$ given by
\[
\l=\eta-\tau= {1\over 2}\, \left( \xi_1-h(s)\,\cos\,\t
\,\sin\Theta , \xi_2 -h(s)\,\sin\t \,\sin\Theta ,
s-h(s)\,\cos\,\Theta \right)\, .
\]
It holds $d\l={h(s)\over 8}\,\sin\Theta\, |s-h(s)\cos\Theta|\,
ds\, d\Theta\, d\t$.

Since $\tau\in \mathcal{B}_{\xi}(\xi + sz)$ we know that
$\taum\leq \left({1\over\sqrt{2}}+{1\over 50}\right)\e$ and hence,
\[
\left(1-\left({1\over\sqrt{2}}+{1\over 50}\right)\right)\e\leq
|\eta-\tau|\leq \e .
\]
Since $(\Theta , \t)\in \mathcal{M}(s , z )$ it holds $|\l|\sim
h(s)$ ($\,|\eta-\tau|\sim\e\,$). Moreover, $\xitau\leq {\e\over
50}$ and the angle $\g$ between $\eta-\xi$ and $\eta-\tau$
satisfies $|\cos\g|\geq C>0$. That is, $|(\eta-\xi)\cdot
(\eta-\tau)|\sim |\eta-\xi|\,|\eta-\tau|$. This says that
$|s-h(s)\cos\Theta|\sim |\l|$, since $z\cdot (\eta-\tau) =
{s-h(s)\cos\Theta\over 2}$. On our domain of integration we have
$|s-h(s)\cos\Theta|\sim s$:
\[
|s-h(s)\cos\Theta|\sim |\l|\sim h(s)\sim s ,
\]
where the condition $h(s)\sim s$ follows from $s\geq 0,9\xim$.

We conclude
\[
G(\xi)\leq C \int_{\mathbb
S^1}\int_{\Rtres}|\q(\l)|^2|\l|^{2\b-1+2\epsi} d\l \, d\s(z)\, .
\]

\qed

%%%%%%%%%%%%%%%%%%%%%%%%%%% Bibliography %%%%%%%%%%%%%%%%%%%%%%%%%%%%%%

\end{document}